%% file: main_HASWME_revisionFinal.tex
\PassOptionsToPackage{unicode}{hyperref}
\PassOptionsToPackage{naturalnames}{hyperref}

\documentclass[a4paper,11pt]{article}

\usepackage{etex}
\usepackage{graphicx}
\usepackage{amsmath}
\usepackage{amssymb}
\usepackage{amsthm}
\usepackage{empheq}
\usepackage{multirow}
\usepackage{makecell}
\usepackage{pict2e}
\usepackage{booktabs}
\usepackage[figuresright]{rotating}
\usepackage{bbm}

\usepackage{wrapfig}
\usepackage{units}
\usepackage{mathtools}  
\usepackage{xfrac}
\usepackage{tikz}
\usepackage[abs]{overpic}
\usepackage{xcolor,varwidth}
\usetikzlibrary{angles,quotes}
\usetikzlibrary{bending}
\tikzset{>=latex} % for LaTeX arrow head
\usepackage{xcolor}
\colorlet{veccol}{green!70!black}
\colorlet{vcol}{green!70!black}
\colorlet{xcol}{blue!85!black}
\colorlet{projcol}{xcol!60}
\colorlet{unitcol}{xcol!60!black!85}
\colorlet{unitcol2}{vcol!60!black!85}
\colorlet{myblue}{blue!70!black}
\colorlet{myred}{red!70!black}
\tikzstyle{vector}=[->,very thick,xcol]
\tikzstyle{mydashed}=[dash pattern=on 2pt off 2pt]

\usepackage{nomencl}   % For nomenclature

\makenomenclature%

\usepackage{textcomp} % nice greek alphabet
\usepackage{pifont}   % Dingbats
\usepackage{booktabs}
\usepackage{enumerate}
\usepackage{physics}
\usepackage{wrapfig}
\usepackage{subcaption}
\usepackage{lineno}
\newcommand\myeq{\stackrel{\mathclap{\normalfont\mbox{\textbf{lin}}}}{=}}
\newtheorem{theorem}{Theorem}[section]

\newtheorem*{theorem*}{Theorem}
\newtheorem*{corollary*}{Corollary}
\newtheorem*{lemma*}{Lemma}
\newtheorem*{hypothesis*}{Hypothesis}
\newtheorem*{proposition*}{Proposition}
\newtheorem*{conjecture*}{Conjecture}

\theoremstyle{definition}

\makeatletter
\newtheorem*{rep@theorem}{\rep@title}
\newcommand{\newreptheorem}[2]{%
\newenvironment{rep#1}[1]{%
 \def\rep@title{#2 \ref{##1}}%
 \begin{rep@theorem}}%
 {\end{rep@theorem}}}
\makeatother

\newreptheorem{theorem}{Theorem}
\newreptheorem{lemma}{Lemma}

\usepackage{algorithm}      % For algorithms
\usepackage{pgfplots}
\pgfplotsset{compat=1.17}
\usepackage{listings}		%For code listings
\usepackage[noend]{algpseudocode}
\usepackage{changepage}     % for changing margins
\usepackage{color}              % for coloring of changed matrix entries
\definecolor{myred}{gray}{0} %{rgb}{1,0,0}
\definecolor{light-gray}{gray}{0.8}
\usepackage{epstopdf}
\usepackage{comment}

\usepackage{graphics} %% add this and next lines if pictures should be in esp format
\usepackage{epsfig} %For pictures: screened artwork should be set up with an 85 or 100 line
\usepackage[ plainpages = false, pdfpagelabels,
	     bookmarks,
	     bookmarksopen = true,
	     bookmarksnumbered = true,
	     breaklinks = true,
	     linktocpage,
	     pagebackref,         % shows hyperrefs for usage of refs
	     colorlinks = false,  % was true
         hypertexnames=false,   %destination with the same identifier error fixed
	     linkcolor = blue,
	     urlcolor  = blue,
	     citecolor = red,
	     anchorcolor = green,
	     hyperindex = true,
	     hyperfigures
	     ]{hyperref}
\usepackage{pdflscape}
\usepackage{overpic}


\usepackage{fancyhdr}
\theoremstyle{remark}
\newtheorem{remark}[theorem]{Remark}

\theoremstyle{definition}

\let\originalleft\left
\let\originalright\right
\renewcommand{\left}{\mathopen{}\mathclose\bgroup\originalleft}
\renewcommand{\right}{\aftergroup\egroup\originalright}

\binoppenalty=\maxdimen % makes line-breaking for inline math formulas impossible
\relpenalty=\maxdimen   % makes line-breaking for inline math formulas impossible

\numberwithin{equation}{section}    % numbering of equations by sections or subsection
\date{\today}             % no date displayed

\title{Capturing vertical information in radially symmetric flow using hyperbolic shallow water moment equations}
\author{
Rik Verbiest\footnote{Corresponding author, email address {r.verbiest@rug.nl}} \footnote{Bernoulli Institute, University of Groningen},
Julian Koellermeier\footnotemark[\value{footnote}]  
}

\begin{document}

\maketitle
%\tableofcontents
%\newpage

\input{axisymmetric/Abstract}
\input{axisymmetric/Introduction}

\input{axisymmetric/HASWME}
\input{axisymmetric/Numerics}
\input{axisymmetric/Conclusion}
\input{axisymmetric/Appendix}

\bibliographystyle{abbrv}
\bibliography{main_HASWME_revisionFinal}

\end{document}

%% file: axisymmetric/Abstract.tex
\begin{abstract}
Models for shallow water flow often assume that the lateral velocity is constant over the water height. The recently derived shallow water moment equations are an extension of these standard shallow water equations. The extended models allow for a vertically changing velocity profile, resulting in more accuracy when the velocity varies considerably over the height of the fluid. Unfortunately, already the one-dimensional models lack global hyperbolicity, an important property of partial differential equations that ensures that disturbances have a finite propagation speed.

In this paper, cylindrical shallow water moment equations are formulated by starting from the cylindrical incompressible Navier-Stokes equations. We formulate two-dimensional axisymmetric Shallow Water Moment Equations by imposing axisymmetry in the cylindrical model. The loss of hyperbolicity is analyzed and a hyperbolic axisymmetric moment model is then derived by modifying the system matrix in analogy to the one-dimensional case, for which the hyperbolicity problem has already been observed and overcome. Numerical simulations with both discontinuous and continuous initial data in a cylindrical domain are performed using a finite volume scheme tailored to the cylindrical mesh. The newly derived hyperbolic model is clearly beneficial as it gives more stable solutions and still converges to the reference solution when increasing the number of moments.

\end{abstract}

%% file: axisymmetric/Introduction.tex
\section{Introduction}
The Shallow Water Equations (SWE) are a set of partial differential equations that describe fluid flows for which the horizontal length scale is much larger than the vertical length scale. Applications can be found in a wide range of scientific fields, such as weather forecasting \cite{weatherForecasting} and free-surface flows like tsunami modelling \cite{tsunamiHakata}. 
%For an example of the SWE applied to simulations of the global atmosphere for weather and climate modelling, more precisely, modelling the flow over a mountain, see \cite{weatherForecasting}. An example of tsunami modelling using the SWE is given in \cite{tsunamiHakata}. 
A crucial feature of the SWE is that the lateral velocity field is constant over the vertical position variable. This is a severe simplification and renders the SWE inaccurate in applications such as dam-break floods \cite{inaccuracy_dambreak} and tsunamis \cite{inaccuracy_tsunami}. Consider, for example, small velocities at the bottom of a river, because of the interaction with the bottom friction. Another example is a large velocity just below the surface due to a strong wind blowing above the free surface.

These applications clearly show that a more flexible system of equations is needed for the modeling of flows with complex motion. For this reason, the so-called \emph{Shallow Water Moment Equations (SWME)} were derived in \cite{SWME}. These equations allow vertical variability in the lateral velocities. The SWME are obtained using the method of moments, which includes an expansion of the lateral velocity in a polynomial basis and a subsequent Galerkin projection to obtain evolution equations for the expansion coefficients. The new system of equations proved to be more accurate than the SWE in numerical simulations \cite{SWME}. The model was recently extended to the non-hydrostatic case in \cite{scholz_dispersion_2023} and generalized to include multi-layer models \cite{Garres-Diaz2023}, similar to \cite{fernandez-nieto2016}.

When modeling flows in rivers or oceans, complex propagation speeds are nonphysical, because waves propagate with real and finite propagation speeds. Unfortunately, already the one-dimensional SWME are not globally hyperbolic \cite{HSWME} leading to propagation speeds with non-zero imaginary part. The loss of hyperbolicity can lead to instabilities in numerical test cases \cite{HSWME, lossOfHyperbolicity}.  The lack of hyperbolicity motivated the derivation of a hyperbolic regularization, the so-called \emph{Hyperbolic Shallow Water Moment Equations (HSWME)} \cite{HSWME}, by modifying the system matrix based on similar approaches from kinetic theory \cite{Fan2016,Koellermeier2020g,Koellermeier2014}, thus guaranteeing global hyperbolicity. The hyperbolicity of the HSWME was proved in one spatial dimension and numerical simulations of the one-dimensional HSWME yielded accurate results \cite{HSWME}. This also allowed for a deeper analysis in \cite{equilibriumStability} and \cite{koellermeier_steady_2022}.

The goal of this paper is to formulate, analyze and simulate moment equations for shallow flows in a cylindrical coordinate system as the necessary step towards a future two-dimensional hyperbolic model. The SWE formulated in cylindrical coordinates are widely used. The reason for this is that for many classical applications such as tsunamis \cite{deb_roy_nonlinear_2007,tsunamiHakata} and tropical cyclones \cite{tropicalCyclones}, a system expressed in cylindrical coordinates is more appropriate. Starting from the cylindrical SWE, an axisymmetric system that models the evolution of flow that only varies in the radial direction can be obtained by setting the angular derivatives to zero. In this way, axisymmetric SWE are classically obtained \cite{axisymmetricGravityCurrents,selfSimilar}. 
We show that the extension to axisymmetric SWME suffers from the same loss of hyperbolicity that is observed in the one-dimensional case \cite{HSWME}. 
As an extension of the one-dimensional case, \emph{Hyperbolic Axisymmetric Shallow Water Moment Equations (HASWME)} are derived for the axisymmetric case in this paper to overcome the hyperbolicity problem. The new axisymmetric models are pseudo-two-dimensional: while there is only a radial dependence of the variables, the model describes two-dimensional flow. In this way, the axisymmetric model can be seen as the next step towards the analysis of a full two-dimensional model.

We therefore prove that the new models have real eigenvalues. Then, we state the hyperbolicity of the new models as proved in \cite{Bauerle_Rotational} for two-dimensional Cartesian models. Next, a finite volume method tailored to a cylindrical grid is described and numerical tests for discontinuous and smooth initial conditions are performed to investigate the accuracy of the new HASWME model. In particular, we show that the hyperbolic regularization results in a numerical solution with less oscillations in a dam break situation. For this test case, it is shown that the hyperbolic model is more accurate for short time horizons. Further, we show that the approximation error decreases when using more moments for smooth initial data, thus indicating convergence. 

The main contributions of this paper are the construction of the axisymmetric SWME and its analysis revealing the loss of hyperbolicity (Section 2), the subsequent derivation of the hyperbolic model called HASWME including an explicit construction of the system matrix for an arbitrary number of moments and an analysis of the eigenvalues of the system matrix (Section 3) and the numerical simulation of cylindrical SWME using a tailored finite volume method (Section 4). 

%A thorough analysis of the stability of the equations also needs to consider the right-hand side source terms and the coupling between these source terms and the transport system matrix. The one-dimensional cartesian HSWME contains both stable and unstable equilibrium \cite{equilibriumStability}, based on a set of structural stability conditions proposed in \cite{stabilityConditions}. The difficulty in an extension of the analysis developed in \cite{equilibriumStability} lies in the forcing terms that are present in the model due to the coordinate transformation. 
%\jk{Equilibrium analysis is not included in the paper any more, so we skip this paragraph}

%% file: axisymmetric/HASWME.tex
\section{Axisymmetric Shallow Water Moment Equations}
\label{ASWME}
This section describes the formulation of axisymmetric shallow water moment equations, as a consistent extension of the Cartesian moment models, derived in \cite{SWME}. Next, we show that the axisymmetric systems lack global hyperbolicity by analysing the eigenvalues of the system matrix.

\subsection{Derivation of reference system}
Formulating the standard incompressible Navier-Stokes equations in cylindrical coordinates %performing an analogous asymptotic analysis as in \cite{SWME}, the following initial system in cylindrical coordinates 
\((r,\theta,z)\in\mathbb{R}^+ \times [0,2\pi] \times [h_b(t,r,\theta),h_s(t,r,\theta)]\), with \(h_b(t,r,\theta)\) and \(h_s(t,r,\theta)\) the time- and position-dependent bottom topography and surface topography, respectively, and assuming a hydrostatic pressure, we obtain 
\begin{align}\label{initialSystemCyl1}
    \partial_r v_r+\frac{v_r}{r}+\frac{1}{r}\partial_{\theta}v_{\theta}+\partial_z w&=0,\\\label{initialSystemCyl2}
    \partial_t v_r+\partial_r(v_r^2)+\frac{1}{r}\partial_{\theta}(v_rv_{\theta})+\partial_z\left( v_r w \right)+\frac{1}{r}\left(v_r^2-v_{\theta}^2\right)&=-\frac{1}{\rho}\partial_r p+\frac{1}{\rho}\partial_z\tau_{rz}+g e_r,\\\label{initialSystemCyl3}
    \partial_t v_{\theta}+\partial_r (v_rv_{\theta})+\frac{1}{r}\partial_{\theta}\left(v_{\theta}^2\right)+\partial_z\left(v_{\theta}w \right)+\frac{2}{r}v_rv_{\theta}&=-\frac{1}{\rho}\frac{1}{r} \partial_{\theta} p+\frac{1}{\rho}\partial_z\tau_{\theta z}+g e_{\theta},\\\label{initialSystemCyl4}
    0&=-\frac{1}{\rho}\partial_z p+g e_z,
\end{align}
where the variables of interest are the water height \(h(t,r,\theta): [0,T] \times \mathbb{R}^+ \times [0,2\pi] \mapsto \mathbb{R} \), the radial velocity \(v_r(t,r,\theta,z): [0,T] \times \mathbb{R}^+ \times [0,2\pi] \times [h_b(t,r,\theta),h_s(t,r,\theta)] \mapsto \mathbb{R}\), and the angular velocity \(v_{\theta}(t,r,\theta,z): [0,T] \times \mathbb{R}^+ \times [0,2\pi] \times [h_b(t,r,\theta),h_s(t,r,\theta)] \mapsto \mathbb{R}\). The deviatoric stress tensor terms in cylindrical coordinates are denoted by \(\tau_{rz}\) and \(\tau_{\theta z}\). From the last equation, we find that the hydrostatic pressure is given by
\begin{equation*}
    p(t,r,\theta,z)=(h_s(t,r,\theta)-z)\rho g e_z.
\end{equation*}

Similar to \cite{SWME}, the variable \(z\) is transformed to \(\xi \in [0,1]\), the so-called $\sigma$-coordinates, using
\begin{equation*}
    \zeta=\frac{z-h_b(t,r,\theta)}{h(t,r,\theta)}.
\end{equation*}
This maps the variables \(v_{r}(t,r,\theta,z)\) and \(v_{\theta}(t,r,\theta,z)\) to \(\Tilde{v}_{r}(t,r,\theta,\zeta)\) and \(\Tilde{v}_{\theta}(t,r,\theta,\zeta)\), respectively, but the tilde will be dropped for readability.
%This effectively maps the vertical variable \(z\) to a variable with a constant domain of \([0,1]\). 
The so-called vertical coupling defined in \cite{SWME} is
\begin{equation*}
    h\omega[h,v_r,v_{\theta}]=-\partial_r\left(h \int_0^{\zeta}v_{r,d}d\Hat{\zeta} \right)-\frac{1}{r}\partial_{\theta}\left(h \int_0^{\zeta}v_{\theta,d}d\Hat{\zeta} \right)- \frac{h}{r}\int_0^{\zeta}v_{r,d}d\Hat{\zeta}, 
\end{equation*}
where \(v_{r,d}=v_{r,d}(t,r,\theta,\zeta)\) and \(v_{\theta,d}=v_{\theta,d}(t,r,\theta,\zeta)\) are the deviations from the average radial velocity defined as \(v_{r,m}(t,r,\theta):=\int^1_0 v_r(t,r,\theta,\zeta) d\zeta\) and the average angular velocity defined as \(v_{\theta,m}(t,r,\theta):=\int^1_0 v_{\theta}(t,r,\theta,\zeta) d\zeta\), respectively.

Using these transformations, the reference system in cylindrical coordinates reads
\footnotesize
\begin{align}\label{refsystemcylindrical1}
    &\partial_t h+\partial_r (hv_{r,m})+\frac{1}{r}\partial_{\theta}(hv_{\theta,m})+\frac{1}{r}hv_{r,m}=0,\\\label{refsystemcylindrical2}
    &\partial_t (hv_r)+\partial_r\left(hv_r^2+\frac{g}{2}e_zh^2\right)+\frac{1}{r}\partial_{\theta}(hv_rv_{\theta})+\partial_{\zeta}\left(h v_r \omega-\frac{1}{\rho}\tau_{rz} \right)+\frac{h}{r}\left(v_r^2-v_{\theta}^2\right)=gh(e_r-e_z\partial_rh_b),\\\label{refsystemcylindrical3}
    &\partial_t (hv_{\theta})+\partial_r (hv_rv_{\theta})+\frac{1}{r}\partial_{\theta}\left(hv_{\theta}^2+\frac{g}{2}e_zh^2\right)+\partial_{\zeta}\left(h v_{\theta}\omega -\frac{1}{\rho}\tau_{\theta z} \right)+\frac{2h}{r}v_rv_{\theta}=gh\left(e_y-\frac{1}{r}e_z\partial_{\theta}h_b\right). 
\end{align}
\normalsize
Note the geometric source terms containing a factor \(\frac{1}{r}\) due to the cylindrical coordinates. These factors make the problem space dependent. %The forcing terms contain a factor \(\frac{1}{r}\), so care has to be taken when doing numerical simulations. 

Furthermore, the system in cylindrical coordinates has some similarity with the one-dimensional system analyzed in \cite{HSWME}, especially if the flow is considered uniform in the angular direction $\theta$ (the axisymmetric case). The hyperbolic regularization later in this paper is therefore based on the findings of the one-dimensional case.

\subsection{Moment equations}\label{section:momentequations}
Analogously to the derivation of the one-dimensional SWME in \cite{SWME}, moment equations can be derived for the reference system \eqref{refsystemcylindrical1}-\eqref{refsystemcylindrical3}. First, the deviation of the radial velocity \(v_{r}(t,r,\theta,\zeta)\) and the angular velocity \(v_{\theta}(t,r,\theta,\zeta)\) from their means \(v_{r,m}(t,r,\theta)\) and \(v_{\theta,m}(t,r,\theta)\), respectively, is modelled as a polynomial expansion:
\begin{align*}
    v_r(t,r,\theta,\zeta)&=v_{r,m}(t,r,\theta)+\sum_{j=1}^{N_r}\alpha_j(t,r,\theta)\phi_j(\zeta),\\
    v_{\theta}(t,r,\theta,\zeta)&=v_{\theta,m}(t,r,\theta)+\sum_{j=1}^{N_{\theta}}\gamma_j(t,r,\theta)\phi_j(\zeta),
\end{align*}
where \(\phi_j(\zeta)\) are the shifted and normalized Legendre polynomials, defined by:
\begin{equation}\label{legendre}
    \phi_j(\zeta)=\frac{1}{j!}\frac{d^j}{d\zeta^j}\left( \zeta - \zeta^2 \right)^j.
\end{equation}
The polynomials \eqref{legendre} form an orthogonal basis on $[0,1]$ satisfying the orthogonality relation
\begin{equation*}
    \int^1_0 \phi_m(\zeta) \phi_n (\zeta) d\zeta = \frac{1}{2n+1}\delta_{mn}.
\end{equation*}
The coefficients \(\alpha_j(t,r,\theta): [0,T] \times \mathbb{R}^+ \times [0,2\pi] \mapsto \mathbb{R} \), with \(j \in [1,\ldots, N_r]\), in radial direction, and \(\gamma_j(t,r,\theta): [0,T] \times \mathbb{R}^+ \times [0,2\pi]\mapsto \mathbb{R}\), with \(j \in [1,\ldots, N_{\theta}]\), in angular direction, are the basis coefficients. They are functions of the time \(t\) and the two-dimensional space \((r,\theta)\). The expansion orders in radial direction and in angular direction are denoted by \(N_r\) and \(N_{\theta}\), respectively. Note that \(N_r\) and \(N_{\theta}\) do not have to be equal. An advantage of this approach is its flexibility; a larger \(N_r\) and/or a larger \(N_{\theta}\) allows for modeling more complex flow, while the choice \(N_r=N_{\theta}=0\) leads to a constant velocity profile, as in the SWE.  

Assuming Newtonian flow, we can close the system with
\begin{equation*}
    \frac{1}{\rho}\tau_{rz}=\frac{\nu}{h} \partial_{\zeta}v_r,\qquad \frac{1}{\rho}\tau_{\theta z}=\frac{\nu}{h}\partial_{\zeta}v_{\theta}
\end{equation*}
and use (weakly enforced) boundary conditions according to \cite{SWME} 
\begin{equation*}
    -\frac{\nu}{h}[\partial_{\zeta} v_r]^{\zeta=1}_{\zeta=0}=\frac{\nu}{\lambda}v_r|_{\zeta = 0}, \qquad-\frac{\nu}{h}[\partial_{\zeta} v_{\theta}]^{\zeta=1}_{\zeta=0}=\frac{\nu}{\lambda}v_{\theta}|_{\zeta = 0},
\end{equation*}
with slip length \(\lambda\) and kinematic viscosity \(\nu\). 

The full cylindrical SWME are derived by multiplying \eqref{refsystemcylindrical2} with the Legendre polynomials \(\phi_j(\zeta)\) (\(j=0,\ldots,N_r\)) and integrating with respect to \(\zeta\), and by multiplying \eqref{refsystemcylindrical3} with the \(\phi_j(\zeta)\) (\(j=0,\ldots,N_{\theta}\)) and integrating with respect to \(\zeta\) \cite{SWME}. The mass balance equation \eqref{refsystemcylindrical1} completes the moment equations. The resulting equations are given below.

\paragraph{Mass balance equation:}
\begin{equation*}
    \partial_{t}h+\frac{1}{r}hv_{r,m}+\partial_r\left( hv_{r,m}\right)+\frac{1}{r}\partial_{\theta}\left(hv_{\theta,m}\right)=0.
\end{equation*}

\paragraph{Projected momentum balance equations:}
\begin{multline}\label{eq:ProjectedMomentumBalanceRadial}
    \partial_t(hv_{r,m})+\partial_r \underbrace{\left( h\left( v_{r,m}^2+\sum_{j=1}^{N_r}\frac{\alpha_j^2}{2j+1}\right) +\frac{g}{2}e_zh^2\right)}_{:=F_r^0}+\frac{1}{r}\partial_{\theta}\underbrace{\left( h\left( v_{r,m}v_{\theta,m}+\sum_{j=1}^{min\{N_r,N_{\theta}\}}\frac{\alpha_j\gamma_j}{2j+1} \right) \right)}_{:=F_{\theta}^0}\\ =\frac{1}{r}\underbrace{h\left( -v_{r,m}^2+v_{\theta,m}^2-\sum_{j=1}^{N_r}\frac{\alpha_j^2}{2j+1}+\sum_{j=1}^{N_{\theta}}\frac{\gamma_j^2}{2j+1} \right)}_{:=G_0}\underbrace{-\frac{\nu}{\lambda}\left( v_{r,m}+\sum_{j=1}^{N_r}\alpha_j \right)}_{:=S_0}+\underbrace{gh(e_r-e_z\partial_r h_b)}_{\text{Bottom topography}},
\end{multline}

\begin{multline}\label{eq:ProjectedMomentumBalanceAngular}
    \partial_t(hv_{\theta,m})+\partial_r\underbrace{\left( h\left( v_{r,m}v_{\theta,m}+\sum_{j=1}^{min\{N_r,N_{\theta}\}}\frac{\alpha_j\gamma_j}{2j+1} \right) \right)}_{:=\Tilde{F}_r^0}+\frac{1}{r}\partial_{\theta} \underbrace{\left( h\left( v_{\theta,m}^2+\sum_{j=1}^{N_{\theta}}\frac{\gamma_j^2}{2j+1}\right) +\frac{g}{2}e_zh^2\right)}_{:=\Tilde{F}_{\theta}^0} \\=\frac{1}{r}\underbrace{\left(-2h\left( v_{r,m}v_{\theta,m}+\sum_{j=1}^{min\{N_r,N_{\theta}\}}\frac{\alpha_j\gamma_j}{2j+1} \right)\right)}_{:=\Tilde{G}_0}
    \underbrace{-\frac{\nu}{\lambda}\left( v_{\theta,m}+\sum_{j=1}^{N_{\theta}}\gamma_j \right)}_{\Tilde{S}_0}+\underbrace{gh\left(e_{\theta}-e_z\frac{1}{r}\partial_{\theta}h_b\right)}_{\text{Bottom topography}}.
\end{multline}

\paragraph{Moment equations:} 
The equations for the coefficients of the radial velocity profile \(\alpha_i, i=1,\ldots, N_r\), are given by:
\begin{equation*}
    \partial_t(h\alpha_i)+\partial_r F_r^i+\frac{1}{r}\partial_{\theta}F_{\theta}^i=Q_r^i\frac{\partial V}{\partial r}+\frac{1}{r}Q_{\theta}^i\frac{\partial V}{\partial \theta}+\frac{1}{r}G_i+S_i,
\end{equation*}
with 
\begin{align*}
    F_r^i&=h\left( 2v_{r,m}\alpha_i+\sum_{j,k=1}^{N_r}A_{ijk}\alpha_j\alpha_k \right),~
    F_{\theta}^i=h\left( v_{r,m}\gamma_i+v_{\theta,m}\alpha_i+\sum_{j=1}^{N_r}\sum_{k=1}^{N_{\theta}}A_{ijk}\alpha_j\gamma_k \right),\\[4pt]
    Q_r^i\frac{\partial V}{\partial r}&=v_{r,m}\partial_r(h\alpha_i)-\sum_{j,k=1}^{N_r}B_{ijk}\alpha_k\partial_r(h\alpha_j),~
    Q_{\theta}^i\frac{\partial V}{\partial \theta}=v_{r,m}\partial_{\theta}(h\gamma_i)-\sum_{j=1}^{N_{\theta}}\sum_{k=1}^{N_r}B_{ijk}\alpha_k \partial_{\theta}(h\gamma_j),\\[4pt]
    G_i&= h\left( -v_{r,m}\alpha_i+2v_{\theta,m}\gamma_i-\sum_{j,k=1}^{N_r}A_{ijk}\alpha_j\alpha_k + \sum_{j,k=1}^{N_{\theta}}A_{ijk}\gamma_j\gamma_k -\sum_{j,k=1}^{N_r}B_{ijk}\alpha_k \alpha_j \right),\\[4pt]
    S_i&=-(2i+1)\frac{\nu}{\lambda}\left( v_{r,m}+\sum_{j=1}^{N_r}\alpha_j\left( 1+\frac{\lambda}{h}C_{ij} \right) \right).
\end{align*}
\(F_r^i\) and \(\frac{1}{r} F_{\theta}^i\) are denoted as the conservative flux, while \(Q_r^i\) and \(\frac{1}{r} Q_{\theta}^i\) contain the non-conservative flux. %Note that \(i\) is a superscript in this notations, to indicate that they belong to the equation for the \(i\)th basis function \(\alpha_i\). 
The geometric source terms due to the formulation in cylindrical coordinates, are grouped in \(\frac{1}{r}G_i\). \(S_i\) is the friction source term, which contains the parameters \(\lambda\) and \(\nu\).
The equations for the coefficients of the angular velocity profile \(\gamma_i, i=1,\ldots,N_{\theta},\) read:
\begin{equation*}
    \partial_t(h\gamma_i)+\partial_r \Tilde{F}_r^i+\frac{1}{r}\partial_{\theta}\Tilde{F}_{\theta}^i=\Tilde{Q}_r^i\frac{\partial V}{\partial r}+\frac{1}{r}\Tilde{Q}_{\theta}^i\frac{\partial V}{\partial \theta}+\frac{1}{r}\Tilde{G_i}+\Tilde{S_i},
\end{equation*}
with 
\begin{align*}
    \Tilde{F}_r^i&=h\left(v_{r,m}\gamma_i+v_{\theta,m}\alpha_i+\sum_{j=1}^{N_r}\sum_{k=1}^{N_{\theta}}A_{ijk}\alpha_j\gamma_k \right),~
    \Tilde{F}_{\theta}^i=h\left( 2v_{\theta,m}\gamma_i+\sum_{j,k=1}^{N_{\theta}}A_{ijk}\gamma_j\gamma_k \right),
\end{align*}

\begin{align*}
    \Tilde{Q}_r^i\frac{\partial V}{\partial r}&=v_{\theta,m}\partial_r(h\alpha_i)-\sum_{j=1}^{N_r}\sum_{k=1}^{N_{\theta}}B_{ijk}\gamma_k\partial_r(h\alpha_j),~
    \Tilde{Q}_{\theta}^i\frac{\partial V}{\partial \theta}=v_{\theta,m}\partial_{\theta}(h\gamma_i)-\sum_{j,k=1}^{N_{\theta}}B_{ijk}\gamma_k\partial_{\theta}(h\gamma_j),\\[4pt]
    \Tilde{G}_i&=-h\left( 2v_{r,m}\gamma_i+v_{\theta,m}\alpha_i+\sum_{j=1}^{N_r}\sum_{k=1}^{N_{\theta}}\left(2A_{ijk}+B_{ijk}\right)\alpha_j\gamma_k \right),\\[4pt]
    \Tilde{S}_i&=-(2i+1)\frac{\nu}{\lambda}\left(v_{\theta,m}+ \sum_{j=1}^{N_{\theta}}\gamma_j\left( 1+\frac{\lambda}{h}C_{ij} \right) \right).
\end{align*}
\(\Tilde{F}_r^i\) and \(\frac{1}{r}\tilde{F}_{\theta}^i\) are denoted as the conservative flux, while \(\Tilde{Q}_r^i\) and \(\frac{1}{r}\Tilde{Q}_{\theta}^i\) contain the non-conservative flux. %Note that \(i\) is a superscript in this notations, to indicate that they belong to the equation for the \(i\)th basis coefficient \(\gamma_i\). 
The geometric source terms due to the formulation in cylindrical coordinates are grouped in \(\frac{1}{r}\Tilde{G}_i\). \(\Tilde{S}_i\) is the friction source term, which contains the parameters \(\lambda\) and \(\nu\).

Neglecting the effects of the bottom topography in the momentum balance equations \eqref{eq:ProjectedMomentumBalanceRadial} and \eqref{eq:ProjectedMomentumBalanceAngular}, the full system can be written in compact form as
\begin{equation}\label{systemFormCylindrical}
    \frac{\partial V}{\partial t}+A_r\frac{\partial V}{\partial r} +A_{\theta}\frac{1}{r}\frac{\partial V}{\partial \theta}=G(V)+S(V),
\end{equation}
with 
\begin{align*}
    &V=(h,hv_{r,m},h\alpha_1,\ldots,h\alpha_{N_r},hv_{\theta,m},h\gamma_1,\ldots,h\gamma_{N_\theta}), \quad A_r=\frac{\partial F_r}{\partial V}-Q_r, \quad A_{\theta}=\frac{\partial F_{\theta}}{\partial V}-Q_{\theta}, \\[5pt]
    &F_r = \left( hv_{r,m}, F_r^0, \ldots, F_r^{N_r}, \Tilde{F}_r^0, \ldots, \tilde{F}_r^{N_{\theta}} \right), \quad F_{\theta} = \left( hv_{\theta,m}, F_\theta^0, \ldots, F_\theta^{N_r}, \Tilde{F}_\theta^0, \ldots, \tilde{F}_\theta^{N_{\theta}} \right), \\[5pt]
    &Q_r = \left( 0, Q_r^0, \ldots, Q_r^{N_r}, \Tilde{Q}_r^0, \ldots, \tilde{Q}_r^{N_{\theta}} \right), \quad Q_{\theta} = \left( 0, Q_\theta^0, \ldots, Q_\theta^{N_r}, \Tilde{Q}_\theta^0, \ldots, \tilde{Q}_\theta^{N_{\theta}} \right), \\[5pt]
    &G(V) = \frac{1}{r}\left(-hv_{r,m}, G_0, \ldots, G_{N_r}, \Tilde{G}_0, \ldots, \Tilde{G}_{N_{\theta}}\right), \quad S(V) = (0, S_0, \ldots, S_{N_r}, \Tilde{S}_0, \ldots, \Tilde{S}_{N_{\theta}}).
\end{align*}
The vector V contains the unknown variables, in convective form. \(A_r\) and \(A_{\theta}\) are the one-dimensional system matrices in radial direction and angular direction, respectively. %Note that \(\frac{1}{r}\) is factored out in \(F_{\theta}\) and \(Q_{\theta}\) to show their dependence on this factor explicitly. 
\(G(V)\) is a vector containing the geometric source terms, while the vector \(S(V)\) contains the friction source terms. This system can also be written in the following form:

\begin{equation}\label{systemFormCylindrical2}
    \frac{\partial V}{\partial t}+\frac{1}{r}\frac{\partial(rF_r)}{\partial r}+Q_r \frac{1}{r}\frac{\partial (rV)}{\partial r}+\frac{1}{r}\frac{\partial F_\theta}{\partial \theta}+Q_\theta \frac{1}{r}\frac{\partial V}{\partial \theta}=\overline{G}(V)+S(V),
\end{equation}
with new geometric source term \(\overline{G}(V)\), given by

\begin{equation*}
   \overline{G}(V)=\frac{h}{r}
   \begin{pmatrix}
       0 \\
        v_{\theta,m}^2+\sum_{j=1}^{N_{\theta}}\frac{\gamma_j^2}{2j+1}  \\
        2v_{\theta,m}\gamma_1+\sum_{j,k=1}^{N_{\theta}}A_{1jk}\gamma_j\gamma_k  \\
       \vdots \\
        2v_{\theta,m}\gamma_{N_r}+\sum_{j,k=1}^{N_{\theta}}A_{N_rjk}\gamma_j\gamma_k  \\
       - v_{r,m}v_{\theta,m}-\sum_{j=1}^{\min\{N_r,N_{\theta}\}} \frac{\alpha_j \gamma_j}{2j+1}  \\
       - v_{r,m}\gamma_1-v_{\theta,m}\alpha_1-\sum_{j=1}^{N_r}\sum_{k=1}^{N_{\theta}} A_{1jk}\alpha_j\gamma_k  \\
       \vdots \\
        -v_{r,m}\gamma_{N_{\theta}}-v_{\theta,m}\alpha_{N_{\theta}}-\sum_{j=1}^{N_r}\sum_{k=1}^{N_{\theta}} A_{N_{\theta}jk}\alpha_j\gamma_k 
   \end{pmatrix} .
\end{equation*}
By using the divergence operator in polar coordinates, i.e.,

\begin{equation}
    \nabla_{\text{polar}} \cdot (w_r,w_\theta):=\frac{1}{r}\frac{\partial (r w_r)}{\partial r}+\frac{1}{r}\frac{\partial w_\theta}{\partial \theta},
\end{equation}
we can write this more compactly as  

\begin{equation}\label{systemFormCylindrical3}
    \frac{\partial V}{\partial t}+\nabla_{\text{polar}} \cdot (F_r,F_\theta)+Q_r \frac{1}{r}\frac{\partial (rV)}{\partial r}+Q_\theta \frac{1}{r}\frac{\partial V}{\partial \theta}=\overline{G}(V)+S(V).
\end{equation}
From \eqref{systemFormCylindrical3} it follows that the main difference between the cylindrical SWME \eqref{systemFormCylindrical} and the Cartesian SWME (\cite{SWME}) are the geometric source terms due to the formulation in a cylindrical coordinate system.  

\subsection{Axisymmetric Shallow Water Moment Equations (ASWME)}
System \eqref{systemFormCylindrical} carries information in both radial and angular direction. In many situations, however, waves in the fluid propagate predominantly only in radial direction, see Figure \ref{fig:axisymmetricFlow}. %Fluid flows that satisfy this description are said to possess axisymmetry. 
A classical example of such axisymmetric cases are tsunamis \cite{tsunamiHakata,tobias_idealized_2011}. In \cite{tropicalCyclones}, axisymmetric shallow water equations are used to understand the intensity of tropical cyclones. Axisymmetric currents have also been thoroughly studied using the axisymmetric shallow water equations, see, e.g., \cite{axisymmetricGravityCurrents} and \cite{selfSimilar}. An axisymmetric moment system can be obtained starting from \eqref{systemFormCylindrical}.

%\begin{wrapfigure}{r}{4cm} \jk{all figures should be normally inserted in the text, not wrapped, as this is very uncommon in scientific papers.}
\begin{figure}[ht]
    \centering
    \resizebox{6cm}{4.75cm}{
        \def\xmax{2.0}
        \def\ul{0.6}
        \def\R{1.7}
        \begin{tikzpicture}
            \def\ang{43}
            \coordinate (O) at (0,0);
            \coordinate (P1) at (0:\R);
            \coordinate (P2) at (45:\R);
            \coordinate (P3) at (90:\R);
            \coordinate (P4) at (135:\R);
            \coordinate (P5) at (180:\R);
            \coordinate (P6) at (225:\R);
            \coordinate (P7) at (270:\R);
            \coordinate (P8) at (315:\R);
            \coordinate (P9) at (22.5:\R);
            \coordinate (P10) at (67.5:\R);
            \coordinate (P11) at (112.5:\R);
            \coordinate (P12) at (202.5:\R);
            \coordinate (P13) at (247.5:\R);
            \coordinate (P14) at (292.5:\R);
            \coordinate (P15) at (157.5:\R);
            \coordinate (P16) at (337.5:\R);
            \coordinate (X) at (\xmax,0);
            \coordinate (R) at (\ang:\R);
            \draw[->,line width=0.9] (-\xmax,0) -- (1.08*\xmax,0) node[right] {$x$};
            \draw[->,line width=0.9] (0,-\xmax) -- (0,1.08*\xmax) node[left] {$y$};
            \draw[vector] (O) -- (P1);
            \draw[vector] (O) -- (P2);
            \draw[vector] (O) -- (P3);
            \draw[vector] (O) -- (P4);
            \draw[vector] (O) -- (P5);
            \draw[vector] (O) -- (P6);
            \draw[vector] (O) -- (P7);
            \draw[vector] (O) -- (P8);
            \draw[vector] (O) -- (P9);
            \draw[vector] (O) -- (P10);
            \draw[vector] (O) -- (P11);
            \draw[vector] (O) -- (P12);
            \draw[vector] (O) -- (P13);
            \draw[vector] (O) -- (P14);
            \draw[vector] (O) -- (P15);
            \draw[vector] (O) -- (P16);
            \draw[color=red!60, very thick](0,0) circle (1.7);
            \draw[color=red!60, thick](0,0) circle (1.2);
            \draw[color=red!60, very thick](0,0) circle (0.7);
        \end{tikzpicture}
    }
   \caption{Axisymmetric flow.}
   \label{fig:axisymmetricFlow}
\end{figure}
%\end{wrapfigure}

In the cylindrical system, axisymmetric flow does not depend on \(\theta\). %standing in the centre of the flow, one will see the same no matter which direction one is looking at. 
Mathematically, all derivatives with respect to \(\theta\) are zero. Note that this does not mean that there is no angular velocity. %The angular velocity just does not depend on \(\theta\), just like the radial velocity and the height. 
This is a strong simplification, but as argued above there are many phenomena that can be modelled with an axisymmetric system. 
One main advantage of this simplified model is its resemblance to the one-dimensional models, see \cite{HSWME} and \cite{SWME}.

In the axisymmetric case, \eqref{systemFormCylindrical} reduces to
\begin{equation}\label{systemFormCylindricalReduced}
    \frac{\partial V}{\partial t}+A_A\frac{\partial V}{\partial r} = G(V) + S(V),
\end{equation}
where \(A_A:=A_r=\frac{\partial F_r}{\partial V}-Q_r\). We call \eqref{systemFormCylindricalReduced} the \emph{Axisymmetric Shallow Water Moment Equations} (ASWME). The axisymmetric system with radial order \(N_r\) and angular order \(N_{\theta}\) is called the \((N_r,N_{\theta})\)th order axisymmetric system. The system matrix of the \((N_r,N_{\theta})\)th order axisymmetric system is denoted by \(A_A^{(N_r,N_\theta)}\). %The subscript \(A\) indicates that \(A_A^{(N_r,N_\theta)}\) is the system matrix belonging to the axisymmetric model. 
Without loss of generality, one particular class of the ASWME is considered in this paper: axisymmetric systems with full velocity expanded, i.e. the \((N,N)\)th order axisymmetric systems. In this case, both the angular velocity and the radial velocity are approximated by a polynomial expansion of the same order \(N\). In the remainder of this paper, we will simply refer to these systems as the \(N\)th order system. The cases of different orders \(N_r\) and \(N_\theta\) are then trivial and left out for brevity. The analytical form of \eqref{systemFormCylindricalReduced} allows for a deeper investigation, where we focus on hyperbolicity first.
%\enlargethispage{3\baselineskip}

\begin{figure}[htb]
    \centering
    \begin{subfigure}{.45\textwidth}
        \includegraphics[width=\textwidth]{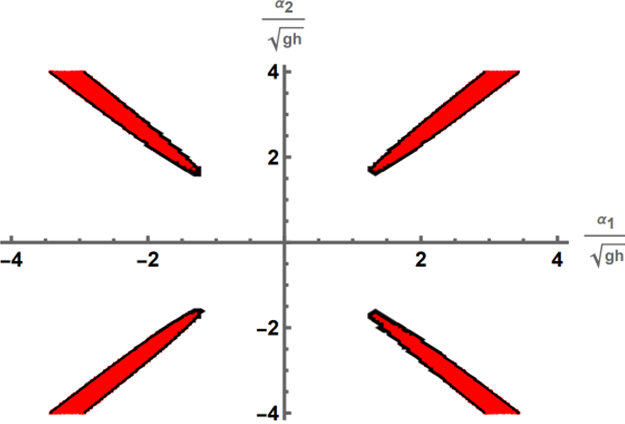}
        \caption{Second order axisymmetric system.}
        \label{fig:22AxisymmetricHypPlot}
    \end{subfigure}
    \hfill
    \begin{subfigure}{.45\textwidth}
        \includegraphics[width=\textwidth]{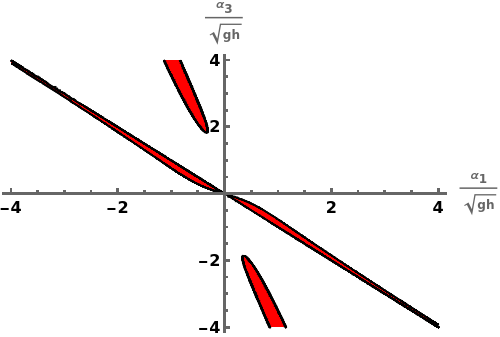}
        \caption{ Third order axisymmetric system.}
        \label{fig:33AxisymmetricHypPlot}
    \end{subfigure}
    \caption{Hyperbolicity region of the second order axisymmetric system (a) and the third order axisymmetric system (b). Red regions indicate loss of hyperbolicity.}
    \label{fig:axisymmetricHypPlots}
\end{figure}

\subsubsection{Hyperbolicity breakdown of axisymmetric system}
In this section, the hyperbolicity of the systems is analyzed. Hyperbolicity is a property of a system of first-order partial differential equations that ensures that information propagates with real and finite propagation speed. It is a requirement for the equations to be robust against small perturbations of the initial data \cite{equilibriumStability,stabilityConditions}. 

In \cite{HSWME}, it is observed that the one-dimensional SWME lack global hyperbolicity, the property of the system matrix to have real and distinct eigenvalues, and an instability that can be related to this loss of hyperbolicity is found in a numerical test case in which shocks are present. 

The \(N\)th order ASWME lack global hyperbolicity, too. This hyperbolicity loss already occurs in the second order system. The system matrix of the second order axisymmetric system is given by:
\begin{equation}
    A_A^{(2,2)}=
    \begin{pmatrix}
        0 & 1 & 0 & 0 & 0 & 0 & 0 \\
        d_1 & 2 v_{r,m} & \frac{2 \alpha _1}{3} & \frac{2 \alpha _2}{5} & 0 & 0 & 0 \\
        d_2 & 2 \alpha _1 & v_{r,m}+\alpha _2 & \frac{3 \alpha _1}{5} & 0 & 0 & 0 \\
        d_3 & 2 \alpha _2 & \frac{\alpha _1}{3} & v_{r,m}+\frac{3 \alpha _2}{7} & 0 & 0 & 0 \\
        \Tilde{d}_1 & v_{\theta ,m} & \frac{\gamma _1}{3} & \frac{\gamma _2}{5} & v_{r,m} & \frac{\alpha _1}{3} & \frac{\alpha _2}{5} \\
        \Tilde{d}_2 & \gamma _1 & \frac{3 \gamma _2}{5} & \frac{\gamma _1}{5} & \alpha _1 & v_{r,m}+\frac{2 \alpha _2}{5} & \frac{2\alpha_1}{5} \\
        \Tilde{d}_3 & \gamma _2 & -\frac{\gamma _1}{3} & \frac{\gamma _2}{7} & \alpha _2 & \frac{2 \alpha _1}{3} & v_{r,m}+\frac{2 \alpha_2}{7}
    \end{pmatrix},
    \label{A22}
\end{equation}
where the entries of the first column are not shown explicitly here for conciseness. They are given in Appendix \ref{app:B}. The eigenvalues of \(A_A^{(2,2)}\) are of the form \(\lambda=v_{r,m}+c\sqrt{gh}\), where \(c\) is any root of the polynomial 
\begin{align*}  \stepcounter{equation}\tag{\theequation}\label{eq:charpol22}
    -\lambda ^7 +
    \frac{74 \alpha _2 \lambda ^6}{35}&+\left(\frac{9 \alpha _1^2}{5}-\frac{177 \alpha _2^2}{245}+1\right) \lambda ^5+\left(-\frac{166 \alpha _2^3}{175}-\frac{234}{175} \alpha _1^2 \alpha _2-\frac{74 \alpha _2}{35}\right) \lambda ^4\\&+\left(-\frac{23 \alpha _1^4}{25}-\frac{486 \alpha _2^2 \alpha _1^2}{1225}-\frac{4 \alpha
   _1^2}{5}+\frac{137 \alpha _2^4}{245}+\frac{324 \alpha _2^2}{245}\right) \lambda ^3 \\ &+\left(\frac{6 \alpha _2^5}{125}+\frac{1908 \alpha _1^2 \alpha _2^3}{6125}-\frac{16 \alpha _2^3}{175}-\frac{6}{35} \alpha _1^4 \alpha _2+\frac{144}{175} \alpha _1^2 \alpha _2\right) \lambda ^2\\&+\left(\frac{3 \alpha _1^6}{25}-\frac{141 \alpha _2^2
   \alpha _1^4}{1225}+\frac{3 \alpha _1^4}{25}+\frac{549 \alpha _2^4 \alpha _1^2}{6125}+\frac{6 \alpha _2^2 \alpha _1^2}{1225}-\frac{323 \alpha _2^6}{6125}-\frac{37 \alpha _2^4}{245}\right) \lambda \\ &+\frac{2 \alpha _2^7}{875}+\frac{54 \alpha _1^2 \alpha _2^5}{6125}+\frac{6 \alpha _2^5}{175}-\frac{278 \alpha _1^4 \alpha
   _2^3}{6125}-\frac{548 \alpha _1^2 \alpha _2^3}{6125}+\frac{6}{175} \alpha _1^6 \alpha _2+\frac{6}{175} \alpha _1^4 \alpha _2.
\end{align*}
Note that, unlike the matrix \(A_A^{(2,2)}\) itself, the characteristic polynomial of this matrix and the eigenvalues do not depend on \(\gamma_1\) and \(\gamma_2\) %As a result, the eigenvalues of \(A_A^{(2,2)}\) do not depend on \(\gamma_1\) and \(\gamma_2\) either. 
due to the lower triangular block structure of \(A_A^{(2,2)}\).
%The mathematical explanation for this is that \(A_A^{(2,2)}\) is a lower triangular block matrix. Its characteristic polynomial is the product of the characteristic polynomial of the diagonal blocks. None of the gammas appear in the diagonal blocks, so the characteristic polynomial does not depend on \(\gamma_1\) and \(\gamma_2\). This is the case in all axisymmetric systems that are considered here, and it can probably be proved that it is also the case for a general axisymmetric system. 

The loss of hyperbolicity is displayed in Figure \ref{fig:22AxisymmetricHypPlot}, which shows the hyperbolicity region, i.e., the region in which \eqref{eq:charpol22} has only real roots, of the second order axisymmetric system. Comparing this hyperbolicty plot to the respective hyperbolicty plots in \cite{HSWME}, we observe that the regions where the axisymmetric system loses hyperbolicity are identical to the 1D Cartesian system. The same is observed for the hyperbolicty plot of the third order axisymmetric system (Figure \ref{fig:33AxisymmetricHypPlot}). This suggests that the lack of hyperbolicity of the axisymmetric systems is closely related to the hyperbolicity loss of the one-dimensional SWME, as we will further analyze in the next section. Loss of hyperbolicity is related to stability problems in the 1D SWME \cite{HSWME} and was also investigated in detail for moment models of rarefied gases \cite{Fan2016,Koellermeier2020g,Koellermeier2014}. It is therefore important to obtain hyperbolic models.

\section{Hyperbolic Axisymmetric Shallow Water Moment Equations}
\label{HASWMESection}
As shown in the previous section, the ASWME clearly lack hyperbolicity. In \cite{HSWME}, the loss of hyperbolicity of the one-dimensional SWME is overcome by modifying the system matrix, based on a similar approach from kinetic theory \cite{Fan2016,Koellermeier2020g,
Koellermeier2014}. We will extend this idea to the quasi one-dimensional ASWME. 

First, note that the first-order axisymmetric system is globally hyperbolic, as can easily be verified by directly computing the system matrix and its eigenvalues. Alternatively, consider the characteristic polynomial \eqref{eq:charpol22} of the second order axisymmetric system with the coefficient \(\alpha_2\) set to zero and observe that \eqref{eq:charpol22} in this case reduces to a polynomial

\begin{equation*}
    -\lambda^7 + \left(\frac{9 \alpha _1^2}{5}+1\right) \lambda ^5+\left(-\frac{23 \alpha _1^4}{25}-\frac{4 \alpha _1^2}{5}\right) \lambda ^3+\left(\frac{3 \alpha _1^6}{25}+\frac{3 \alpha _1^4}{25}\right) \lambda 
\end{equation*}
with real and distinct roots

\begin{equation*}
    \lambda_1=0, \quad \lambda_{2,3} = \pm \sqrt{1+\alpha_1^2}, \quad \lambda_{4,5} = \pm \sqrt{\frac{1}{5}}\alpha_1, \quad \lambda_{6,7} = \pm \sqrt{\frac{3}{5}}\alpha_1,
\end{equation*}
clearly resulting in five distinct, real eigenvalues for the first order system.

The \emph{Hyperbolic Axisymmetric Shallow Water Moment Equations} (HASWME) are derived by linearizing the \(N\)th order system matrix \(A_{A}^{(N,N)}\) around linear deviations from the constant equilibrium velocity profile,  %\((h,v_{r,m},\alpha_1,0,\ldots,0,v_{\theta,m},\gamma_1,0,\ldots,0)\), 
i.e.,
\begin{equation*}
    (h,v_{r,m},\alpha_1,\ldots,\alpha_{N_r},v_{\theta,m},\gamma_1,\ldots,\gamma_{N_{\theta}})\longrightarrow (h,v_{r,m},\alpha_1,0,\ldots,0,v_{\theta,m},\gamma_1,0,\ldots,0).
\end{equation*}
Practically, the higher order coefficients \(\alpha_i\) and \(\gamma_i\), with \(i \geq 2 \), are set to zero in the system matrix. This way, the HASWME system with modified system matrix \(A_{HA}\) is obtained:
\begin{equation}\label{systemFormHyperbolicAxisymmetric}
    \frac{\partial V}{\partial t}+A_{HA}\frac{\partial V}{\partial r} = G(V) + S(V).
\end{equation}
The analytical forms of the hyperbolic axisymmetric systems with full velocity expanded can be derived, allowing for deeper mathematical analysis of the HASWME model.

\begin{remark}
    We note that similar to \cite{koellermeier_steady_2022} and \cite{HSWME}, different ways to perform a hyperbolic regularization exist. Additional regularization terms could be used to construct models with specific eigenvalues. However, we leave this option for potential future work and focus on the standard hyperbolic regularization method here.
\end{remark}

\begin{remark}
    The hyperbolic regularization can be applied to general \((N_r,N_{\theta})\)th order axisymmetric systems by a straightforward generalisation of the construction given in the next section. 
\end{remark}

\subsection{Analytical form of the axisymmetric system with full velocity expanded}
Deriving the full system matrix and then applying the hyperbolic regularization mentioned above, the hyperbolic system matrix of a general axisymmetric system with radial velocity expanded is obtained. We follow the derivation of the 1D system \cite{HSWME} and generalize the procedure for the extended coordinate system to a block matrix structure.

\begin{theorem}\label{NNOrderTheoremMatrixForm}
The HASWME system matrix \(A_{HA}^{(N,N)} \in \mathbb{R}^{(2N+3)\times (2N+3)}\) is given by

\begin{equation}\label{22systemMatrix}
    A_{HA}^{(N,N)}=
    \begin{pmatrix}
        \boldsymbol{A}^{(N,N)} & \boldsymbol{0}^{(N,N)} \\
        \boldsymbol{B}^{(N,N)} & \boldsymbol{C}^{(N,N)}
    \end{pmatrix},
\end{equation}
where
\begin{equation}\label{ANN}
    \boldsymbol{A}^{(N,N)}=
    \begin{pmatrix}
        & 1 & & & &  \\[6pt]
        gh-v_{r,m}^2-\frac{1}{3}\alpha_1^2 & 2v_{r,m} & \frac{2}{3}\alpha_1 & & & \\[6pt]
        -2 v_{r,m}\alpha_1 & 2\alpha_1 & v_{r,m} & \frac{3}{5}\alpha_1 & & & \\[6pt]
        -\frac{2}{3}\alpha_1^2 & & \frac{1}{3}\alpha_1 & v_{r,m} & \ddots & & \\[6pt]
        & & & \ddots & \ddots & \frac{N+1}{2 N+1}\alpha_1 \\[6pt]
        & & & & \frac{N-1}{2N-1}\alpha_1 & v_{r,m}
    \end{pmatrix},
\end{equation}

\begin{equation}\label{BNN}
    \boldsymbol{B}^{(N,N)}=
    \begin{pmatrix}
        -v_{r,m}v_{\theta,m} -\frac{\alpha_1\gamma_1}{3} & v_{\theta,m} & \frac{\gamma_1}{3} & & & \\[6pt]
        -v_{r,m}\gamma_1-v_{\theta,m}\alpha_1 & \gamma_1 & 0  & \frac{2}{5}\gamma_1 & & \\[6pt]
        -\frac{2}{3}\alpha_1\gamma_1 & & -\frac{1}{3}\gamma_1 & 0 & \ddots & \\[6pt]
        & & & \ddots & \ddots & \frac{N}{2N+1}\gamma_1 \\[6pt]
        & & & & -\frac{1}{2N-1}\gamma_1 & 0
    \end{pmatrix},
\end{equation}

\begin{equation}\label{CNN}
    \boldsymbol{C}^{(N,N)}=
    \begin{pmatrix}
        v_{r,m} & \frac{\alpha_1}{3} & & & & \\[6pt]
        \alpha_1 & v_{r,m} & \frac{2}{5}\alpha_1 & & & \\[6pt]
        & \frac{2}{3}\alpha_1 & v_{r,m} & \ddots & & \\[6pt]
        & &  \ddots & \ddots & \frac{N}{2N+1}\alpha_1 \\[6pt]
        & & &  \frac{N}{2N-1}\alpha_1 & v_{r,m}
    \end{pmatrix},
\end{equation}
with \(\boldsymbol{A}^{(N,N)} \in \mathbb{R}^{(N+2)\times (N+2)}\), \(\boldsymbol{B}^{(N,N)}\in \mathbb{R}^{(N+1)\times (N+2)}\) and \(\boldsymbol{C}^{(N,N)} \in \mathbb{R}^{(N+1)\times (N+1)}\) and where all other entries are zero. \(\boldsymbol{0}^{(N,N)} \in \mathbb{R}^{(N+2)\times (N+1)}\) is a zero matrix.
\end{theorem}

\begin{proof}
The proof is an extension of the 1D result in \cite{HSWME} and can be found in Appendix \ref{ch:hyperbolicityproof(N,N)}.
\end{proof}
Note that \(\boldsymbol{A}^{(N,N)}\) is precisely the system matrix of the \(N\)th order 1D HSWME \cite{HSWME}. In \cite{Bauerle_Rotational}, it is claimed that the matrix \(A_{HA}^{(N,N)}\) is diagonalizable with real eigenvalues. However, that proof only holds if \(\alpha_1 \neq 0\), since otherwise there is a repeated eigenvalue and the matrix cannot be diagonalized any more. We will now derive the eigenvalues of \(A_{HA}^{(N,N)}\) for the general case without assuming values of \(\alpha_1\). In particular, we make use of  the fact that the system matrix is a block lower triangular matrix, simplifying the computation of the characteristic polynomial and the eigenvalues below.  

\begin{theorem}\label{NNOrderCharPol_Verbiest}
The HASWME system matrix \(A_{HA}^{(N,N)}\in \mathbb{R}^{(N+3)\times (N+3)}\) has the following characteristic polynomial:
\begin{equation}\label{NNcharpol_Verbiest}
    \chi_{A_{HA}^{(N,N)}}(\lambda)=\left( \left( \lambda -v_{r,m})^2-gh-\alpha_1^2 \right) \right)\cdot \chi_{A_2^{(N,N)}}(\lambda -v_{r,m})\cdot\chi_{A_3^{(N,N)}}(\lambda-v_{r,m}),
\end{equation}
where \(A_2^{(N,N)}\in \mathbb{R}^{N\times N}\) and \(A_3^{(N,N)}\in \mathbb{R}^{(N+1)\times (N+1)}\) are defined as
\begin{equation*}
    A_2^{(N,N)}=
    \begin{pmatrix}
        &  c_1 & &  \\[6pt]
        a_1 & & \ddots &  \\[6pt]
        & \ddots & & c_{N-1} \\[6pt]
        & & a_{N-1} &
    \end{pmatrix}, \qquad
    A_3^{(N,N)}=
    \begin{pmatrix}
         & g_1 & &  \\[6pt]
        f_1 &  & g_2 & & \\[6pt]
        & f_2 &  & \ddots &  \\[6pt]
        & & \ddots &  & g_N \\[6pt]
        & & & f_{N} & 
    \end{pmatrix}
\end{equation*}
with
\begin{align*}
    a_i&=\frac{i}{2i+1}\alpha_1, \quad
    c_i=\frac{i+2}{2i+3}\alpha_1, \qquad i=1,\ldots,N-1,\\[5pt]
    f_i&=\frac{i}{2i-1}\alpha_1,\quad g_i=\frac{i}{2i+1}\alpha_1, \qquad i=1,\ldots,N.
\end{align*}
\end{theorem}

\begin{proof}
The characteristic polynomial of the modified system matrix \(A_H^{(N,N)}\) is by definition 
\begin{equation*}
    \chi_{A_{HA}^{(N,N)}}(\lambda)=\det(A_{HA}^{(N,N)}-\lambda I_{2N+3}) :=\lvert A_{HA}^{(N,N)}-\lambda I_{2N+3} \rvert.
\end{equation*}
Recall that the system matrix has the following form:
\begin{equation*}
    A_{HA}^{(N,N)}=
    \begin{pmatrix}
        \boldsymbol{A}^{(N,N)} & \boldsymbol{0} \\
        \boldsymbol{B}^{(N,N)} & \boldsymbol{C}^{(N,N)}
    \end{pmatrix},
\end{equation*}
where the explicit form of the blocks \(\boldsymbol{A}^{(N,N)},\boldsymbol{B}^{(N,N)}\) and \(\boldsymbol{C}^{(N,N)}\) is given in Theorem \ref{NNOrderTheoremMatrixForm}. This is a lower triangular block matrix. The determinant of a triangular block matrix is given by the determinant of its diagonal blocks, so we have:
\begin{equation*}
\lvert A_{HA}^{(N,N)}-\lambda I_{2N+3} \rvert =\lvert \boldsymbol{A}^{(N,N)}-\lambda I_{N+2} \rvert \cdot \lvert \boldsymbol{C}^{(N,N)}-\lambda I_{N+1} \rvert.
\end{equation*}
According to \cite{HSWME}, the first factor yields:
\begin{equation*}
    \lvert \boldsymbol{A}^{(N,N)}-\lambda I \rvert=\left( \left( \lambda -v_{r,m})^2-gh-\alpha_1^2 \right) \right)\cdot \chi_{A_2^{(N,N)}}(\lambda -v_{r,m}),
\end{equation*}
as \(\boldsymbol{A}^{(N,N)}\) is the system matrix of the \(N\)th order 1D HSWME. The characteristic polynomial of \(\boldsymbol{C}^{(N,N)}\) is analogously given by

\begin{align*}
    \chi_{\boldsymbol{C}^{(N,N)}}(\lambda)=\left| \boldsymbol{C}^{(N,N)} - \lambda I_{N+1} \right| 
    &= \left| A_3^{(N,N)} - (\lambda - v_{r,m})  I_{N+1} \right| \\[5pt] &= \chi_{A_3^{(N,N)}}(\lambda-v_{r,m}).
\end{align*}
This completes the proof.

\end{proof}

It is now proved that the eigenvalues of the characteristic polynomial \eqref{NNcharpol_Verbiest} are real.

\begin{theorem}\label{theorem:eigenvaluesNN_Verbiest}
The eigenvalues of the modified axisymmetric hyperbolic system matrix \(A_{HA}^{(N,N)} \in \mathbb{R}^{(2N+3)\times (2N+3)}\) are the real numbers

\begin{align}
    \lambda_{1,2}&=v_{r,m}\pm \sqrt{gh+\alpha_1^2},\\[5pt]
    \label{eigenvalues1_Verbiest}\lambda_{i+2}&=v_{r,m}+b_i\alpha_1, \qquad i=1,\ldots,N,\\[5pt]
    \label{eigenvalues2_Verbiest}\lambda_{i+2+N}&=v_{r,m}+s_i\alpha_1,\qquad i=1,\ldots,N+1,
\end{align}
with \(b_i \cdot \alpha_1\) the real roots of \(A_2^{(N,N)}\), and \(s_i \cdot \alpha_1\) the real roots of \(A_3^{(N,N)}\), from Theorem \ref{NNOrderCharPol_Verbiest} and where all the \(b_i\)'s are pairwise distinct and where all the \(s_i\)'s are pairwise distinct.
\end{theorem}

\begin{proof}
Recall that the characteristic polynomial of \(A_{HA}^{(N,N)}\) is given by
\begin{equation*}
    \chi_{A_{HA}^{(N,N)}}(\lambda)=\left( \left( \lambda -v_{r,m})^2-gh-\alpha_1^2 \right) \right)\cdot \chi_{A_2^{(N,N)}}(\lambda -v_{r,m})\cdot \chi_{A_3^{(N,N)}}(\lambda-v_{r,m}).
\end{equation*}
From the first factor, we obtain
\begin{equation*}
    \lambda_{1,2}=v_{r,m}\pm \sqrt{gh+\alpha_1^2}.
\end{equation*}
In \cite{SWME}, it is shown that the roots of \(\chi_{A_2^{(N,N)}}(\lambda-v_{r,m})\) have the form \(\lambda_i=v_{r,m}+b_i\cdot \alpha_1\). Moreover, it is proved in \cite{equilibriumStability} that the roots are real.

It remains to prove that the eigenvalues of the matrix \(A_3^{(N,N)}\) are real and distinct when \(\alpha_1\neq 0\). This will be done in an analogous manner as in the proof of the statement that the roots of \(\chi_{A_2^{(N,N)}}(\lambda-v_{r,m})\) are real and distinct when \(\alpha_1 \neq 0\) in \cite{equilibriumStability}. Let \(p_{N+1}(z):=\chi_{A_3^{(N,N)}}\) be the characteristic polynomial of \(A_3^{(N,N)}\). The polynomials \(p_{k}(z)\) follow the recurrence formula
\begin{equation*}
    p_k(z)=-zp_{k-1}(z)-f_{i-1}g_{i-1}p_{k-2}(z), \qquad \text{for any}\quad 2 \leq k \leq N+1,
\end{equation*}
and where \(p_0(z):=1\) and \(p_1(z)=-\lambda\). Then we need to prove that the roots of \(p_{N+1}(z)\) are real and distinct when \(\alpha_1 \neq 0\). To prove this statement, the following propositions are proved by induction:
\begin{enumerate}
    \item \(p_{k+1}(z)\) has \(k+1\) different real roots \(z_1^{(k+1)}<z_2^{(k+1)}<\cdots < z_{k+1}^{(k+1)}\) and the \(k\) roots of \(p_k(z)\) lie between them:
    \begin{equation*}
        z_1^{(k+1)}<z_1^{(k)}<z_2^{(k+1)}< \cdots < z_k^{(k+1)} < z_k^{(k)} < z_{k+1}^{(k+1)}.
    \end{equation*}
    \item If \(k\) is even, the signs of the sequence \(p_k\left( z_1^{(k+1)}\right),p_k\left(z_2^{(k+1)} \right),\ldots, p_k\left(z_{k+1}^{(k+1)} \right)\) are \((+,-,+,-,\ldots,+,-,+)\); otherwise the signs are \((+,-,+,-,\ldots,+,-)\).
\end{enumerate}
The two properties can be easily verified for \(p_0(z)=1\), \(p_1(z)=-z\), \(p_2(z)=z^2-\frac{\alpha_1^2}{3}\) and \(p_3(z)=-z^3+\frac{3\alpha_1^2}{5}z\). This proves the base case of the induction. The induction step is proved in \cite{equilibriumStability}.

\end{proof}

It was rightfully noted in \cite{Bauerle_Rotational} that Theorem \ref{theorem:eigenvaluesNN_Verbiest} is not sufficient to prove hyperbolicity as it is only proved that the eigenvalues \eqref{eigenvalues1_Verbiest}  are distinct and that the eigenvalues \eqref{eigenvalues2_Verbiest} are distinct and not that all the eigenvalues \eqref{eigenvalues1_Verbiest} and \eqref{eigenvalues2_Verbiest} are distinct. This motivated the authors of \cite{Bauerle_Rotational} to prove that the eigenvalues are indeed all distinct if \(\alpha_1 \neq 0\) by explicitly constructing the characteristic polynomial using the Legendre polynomials in a slightly different way as is done in the proof of Theorem \ref{NNOrderCharPol_Verbiest} here. We will now state the results proved in \cite{Bauerle_Rotational} while imposing the condition \(\alpha_1 \neq 0\). 

\begin{theorem}{\cite{Bauerle_Rotational}}\label{22OrderCharPol}
If \(\alpha_1\neq 0\), the HASWME system matrix \(A_{HA}^{(N,N)}\in \mathbb{R}^{(2N+3)\times (2N+3)}\) is real diagonalizable. Its characteristic polynomial is given by:

\begin{equation}\label{NNcharpol_Bauerle}
    \chi_{A_{HA}^{(N,N)}}(\lambda)=\frac{N!(N+1)!}{((2N+1)!!)^2}\alpha_1^{2N+1}P'_{N+1}\left(\frac{\lambda-v_{r,m}}{\alpha_1}\right)P_{N+1}\left(\frac{\lambda -v_{r,m}}{\alpha_1}\right)\left( (\lambda -v_{r,m})^2-gh-\alpha_1^2 \right).
\end{equation}
    Moreover, the eigenvalues are given by

\vspace{-0.5cm}

    \begin{align}
        \lambda_{1,2}&=v_{r,m}\pm \sqrt{gh+\alpha_1^2},\\[5pt]
        \label{eigenvalues1}\lambda_{i+2}&=v_{r,m}+b_i\alpha_1, \qquad i=1,\ldots,N,\\[5pt]
        \label{eigenvalues2}\lambda_{i+2+N}&=v_{r,m}+s_i\alpha_1,\qquad i=1,\ldots,N+1,
    \end{align}
    where \(b_i\) for \(i=1,\ldots,N\) are the roots of the derivative of Legendre polynomial \(P'_{N+1}(\xi)\) and \(s_i\) for \(i=1,\ldots,N+1\) are the roots of Legendre polynomial \(P_{N+1}(\xi)\).

    Thus, if \(\alpha_1\neq0\), then the \(N\)th order HASWME are globally hyperbolic.

\end{theorem}

\begin{proof}
    Proved in \cite{Bauerle_Rotational}.
\end{proof}

\begin{remark}\label{remark:multiplicity_NNHASWME}
According to Theorem \ref{22OrderCharPol}, the roots of the system matrix of a hyperbolic axisymmetric system are real. However, for vanishing first coefficient \(\alpha_1\) all eigenvalues \eqref{eigenvalues1} and \eqref{eigenvalues2} are collapsing to \(v_{r,m}\), which can lead to problems diagonalizing the matrix.
\end{remark}

\begin{remark}
The construction of the \(N\)th order systems with full velocity expanded allow for a generalization to arbitrary systems with order of expansion \(N_r\) in radial direction and order of expansion \(N_{\theta}\) in angular direction, with \(N_r \neq N_{\theta} \neq 0\). A straightforward extension of the theory is that the system matrix will always be a lower triangular matrix, facilitating the computation of the characteristic 
polynomial considerably. Denoting the blocks of the system matrix by
\begin{equation*}
A_{sys}=
    \begin{pmatrix}
        \mathbf{A} & \mathbf{0} \\[5pt]
        \mathbf{B} & \mathbf{C}
    \end{pmatrix},
\end{equation*}
the matrix \(\mathbf{A}\) corresponds to the system matrix in the one-dimensional HSWME \cite{HSWME}. Moreover, the dimension of this matrix is only determined by the radial order \(N_r\). Since the matrix \(A_{sys}\) is lower triangular, block \(\mathbf{B}\), which contains the derivatives with respect to \(h, v_{r,m}\) and \(\alpha_i\), with \(i \in [1,\ldots,N_r]\), that appear in the angular momentum balance equation and the equations for \(\gamma_i\), with \(i \in [1,\ldots,N_{\theta}]\), does not appear in the calculation of the characteristic polynomial. Block \(\mathbf{C}\) contains the derivatives with respect to \( v_{\theta,m}\) and \(\gamma_i\), with \(i \in [1,\ldots,N_{\theta}]\), that appear in the angular momentum balance equation and the equations for \(\gamma_i\), with \(i \in [1,\ldots,N_{\theta}]\). Therefore, the dimension of block \(\mathbf{C}\) only depends on the angular order \(N_{\theta}\). It follows that matrix \(\mathbf{C}\) is precisely the matrix \(C^{(N_\theta,N_\theta)}\) defined in Theorem 
\ref{NNOrderTheoremMatrixForm}. In particular, this means that the hyperbolicity theorems and proofs given in this section for the systems with the same order in both radial and angular direction can be easily generalized to systems with order \(N_r\) in radial direction and order \(N_{\theta}\) in angular direction. 
\end{remark}

% TODO: TRY TO GIVE AN EXPLANANTION

\begin{remark}
The stability properties of the HASWME \eqref{systemFormHyperbolicAxisymmetric} are not solely determined by the system matrix representing the transport part, but also by the dissipation part of the model, which is inscribed in the right-hand side source terms \cite{equilibriumStability,stabilityConditions}. %Hence, it is desirable to also characterize this dissipation behavior in the hyperbolic axisymmetric models. 
An equilibrium stability analysis of the 1D HSWME is performed in \cite{equilibriumStability}, in which equilibrium manifolds of the one-dimensional models are derived and in which a set of stability conditions, proposed in \cite{stabilityConditions}, is verified for each of these manifolds. 
%Due to the similarity between the one-dimensional HSWME and the HASWME, it might seem straightforward to extend the one-dimensional Cartesian analysis to the axisymmetric models. 
For the HASWME, the geometric source terms denoted by \(G(V)\) in \eqref{systemFormHyperbolicAxisymmetric} pose mathematical difficulties for the analytical derivation of explicit expressions for the equilibrium manifolds. The equilibrium analysis is therefore left for future work.
\end{remark}

%% file: axisymmetric/Numerics.tex
\section{Numerical Simulation}
\label{Sec:numerics}

In this section, we first describe the construction of the finite volume scheme that was used to perform the moment model simulations. This includes a description of the geometry and the boundary conditions.

Next, we show that the HASWME overcome the stability problems of the ASWME by considering a cylindrical dam break situation. Then, we demonstrate that the ASWME and the HASWME yield accurate results by considering a test case with smooth initial conditions. Moreover, we observe that, in this test case, the error reduces with respect to increasing order \(N\) throughout all models. The simulations are performed using the axisymmetric systems with full velocity expanded up to order \(N=4\). Simulations with axisymmetric moment models with different orders are left to future work. 

Numerical simulations were performed with the ASWME and the HASWME. An important feature of the cylindrical coordinate system that has to be dealt with is the \(\frac{1}{r}\) factor in some of the terms, leading to a position-dependent flux function, see Equation \eqref{systemFormCylindrical3}, for example. 

\subsection{Finite volume method}

Although the numerical scheme will be used to simulate the 1D ASWME \eqref{systemFormCylindricalReduced} and the 1D HASWME \eqref{systemFormHyperbolicAxisymmetric}, it is necessary to start from a 2D formulation to describe the 2D geometry. A rectangular grid in the \((r,\theta)\)-plane is used. The radial spacing is denoted by \(\Delta r\), while the angular spacing is denoted by \(\Delta \theta\). Both \(\Delta r\) and \(\Delta \theta\) are constant. Let \(I_r\) and \(I_{\theta}\) be the extension of the grid in radial direction and angular direction, respectively. The centers of the control volumes are then given by \((r_i,\theta_j)\), with \(r_i=r_1+i\frac{\Delta r}{I_r}\) and \(\theta_j=\theta_1+j\frac{\Delta \theta}{\theta_r}\), for \(i=1,\ldots,I_r\) and \(j=1,\ldots,I_\theta\). The geometry of a cell \(V_{ij}\) in the mesh is illustrated in Figure \ref{fig:GridCell}. The positions of the faces of a cell are denoted by half indices: the inward face of the cell \((r_i,\theta_i)\) in Figure \ref{fig:GridCell} along boundary \(b_{i-\frac{1}{2},j}\) has radial coordinate \(r_{i-\frac{1}{2}}:=r_i-\frac{\Delta r}{2}\), for example. 
\begin{figure}[h!]
    \centering
    \begin{subfigure}[b]{0.46\textwidth}
        \centering
    \begin{tikzpicture}[>=latex]
        \clip (1.2,5) rectangle + (6.,-6.);
        \coordinate (1) at (0:3);
        \coordinate (2) at (0:5);
        \coordinate (3) at (90:3);
        \coordinate (4) at (90:5);
        \coordinate (5) at (1,0);
        \coordinate (6) at (1,-1);

        \draw[gray] (0,0) -- (0:6);
        \draw[gray] (0,0) -- (40:6);
        \draw[gray] (3,0) arc (0:60:3);
        \draw[gray] (5,0) arc (0:60:5);
        \draw[gray] (3,0) arc (0:-15:3);
        \draw[gray] (5,0) arc (0:-15:5);

        \draw[very thick] (3,0) arc (0:40:3);
        \draw[very thick] (5,0) arc (0:40:5);
        \draw[very thick] (3,0) -- (5,0);
        \draw[very thick] (2.29813332936
        ,1.92836282906) -- (3.83022221559,3.21393804843);
        \node[] at (4,-0.5) {\(b_{i,j-\frac{1}{2}}\)};
        \node[] at (2.7,3) {\(b_{i,j+\frac{1}{2}}\)};
        \node[] at (2.2,1) {\(b_{i-\frac{1}{2},j}\)};
        \node[] at (5.3,1.5) {\(b_{i+\frac{1}{2},j}\)};

        \node[] at (3.5,1.5) {\(V_{i,j}\)};
    \end{tikzpicture}
    \caption{Grid cell in cylindrical mesh.}
    \label{fig:GridCell}
    \end{subfigure}
    \hfill
    \begin{subfigure}[b]{0.46\textwidth}
    \centering
        \resizebox{4cm}{4cm}{
\begin{tikzpicture}
    \coordinate (O) at (0,0);
    \draw (O) circle (3);
    \draw (O) circle (1);
    \draw[line width=0.5mm,color=black] (-1,0) -- (1,0);
    \draw[line width=0.5mm,color=black] (-0.96,0.2) -- (0.96,0.2);
    \draw[line width=0.5mm,color=black] (-0.96,-0.2) -- (0.96,-0.2);
    \draw[line width=0.5mm,color=black] (-0.9,0.4) -- (0.9,0.4);
    \draw[line width=0.5mm,color=black] (-0.9,-0.4) -- (0.9,-0.4);
    \draw[line width=0.5mm,color=black] (-0.8,0.6) -- (0.8,0.6);
    \draw[line width=0.5mm,color=black] (-0.8,-0.6) -- (0.8,-0.6);
    \draw[line width=0.5mm,color=black] (-0.6,0.8) -- (0.6,0.8);
    \draw[line width=0.5mm,color=black] (-0.6,-0.8) -- (0.6,-0.8);

    \coordinate (A1) at (2.5,0); \coordinate (A2) at (4,-1);
    \coordinate (B1) at (1.77,1.77); \coordinate (B2) at (3.5,2);
    \coordinate (C1) at (0,2.5); \coordinate (C2) at (1,4);
    \coordinate (D1) at (-1.77,1.77); \coordinate (D2) at (-3,3.25);
    \coordinate (E1) at (-2.5,0); \coordinate (E2) at (-4,1);
    \coordinate (F1) at (-1.77,-1.77); \coordinate (F2) at (-3.25,-2.75);
    \coordinate (G1) at (0,-2.5); \coordinate (G2) at (-1,-4);
    \coordinate (H1) at (1.77,-1.77); \coordinate (H2) at (2.5,-3.55);
    
    \draw [-{>[sep=2pt]},line width=0.5mm] (A1) to [bend left=45] (A2);
    \draw [-{>[sep=2pt]},line width=0.5mm] (B1) to [bend left=45] (B2);
    \draw [-{>[sep=2pt]},line width=0.5mm] (C1) to [bend left=45] (C2);
    \draw [-{>[sep=2pt]},line width=0.5mm] (D1) to [bend left=45] (D2);
    \draw [-{>[sep=2pt]},line width=0.5mm] (E1) to [bend left=45] (E2);
    \draw [-{>[sep=2pt]},line width=0.5mm] (F1) to [bend left=45] (F2);
    \draw [-{>[sep=2pt]},line width=0.5mm] (G1) to [bend left=45] (G2);
    \draw [-{>[sep=2pt]},line width=0.5mm] (H1) to [bend left=45] (H2);
    \end{tikzpicture}
}
\caption{Axisymmetric domain with solid wall and outflow boundary conditions.}
    \label{fig:boundaryIllustration}
\end{subfigure}
    \caption{Graphical representation of a grid cell (a) and the boundary conditions (b).}
    \label{fig:gridAndBoundary}
\end{figure}
The piecewise boundaries of the cell \(V_{ij}\) are denoted by \(b_{i,j-\frac{1}{2}}\), \(b_{i+\frac{1}{2},j}\), \(b_{i,j+\frac{1}{2}}\) and \(b_{i-\frac{1}{2},j}\). We denote the respective lengths by \(l_{i,j-\frac{1}{2}}\), \(l_{i+\frac{1}{2},j}\), \(l_{i,j+\frac{1}{2}}\) and \(l_{i-\frac{1}{2},j}\), and the volume of cell \(V_{ij}\) is denoted by \(\left| V_{ij} \right|\). 
Following \cite{LeVeque_2002}, we use a finite volume method of the form: 
\begin{align}
    V_{ij}^{n+1}=V_{ij}^n&-\frac{\Delta t}{|V_{ij}|}\left( l_{i-\frac{1}{2},j} \mathcal{A}^+\Delta V_{i-\frac{1}{2},j}+l_{i+\frac{1}{2},j}\mathcal{A}^-\Delta V_{i+\frac{1}{2},j}\right) \notag\\ \label{eq:finiteVolumeScheme}
    &- \frac{\Delta t}{|V_{ij}|}\left(l_{i,j-\frac{1}{2}} \mathcal{B}^+\Delta V_{i,j-\frac{1}{2}}+l_{i,j+\frac{1}{2}}\mathcal{B}^-\Delta V_{i,j+\frac{1}{2}} \right),
\end{align}
where \(\mathcal{A}^{\pm}\Delta V_{i-\frac{1}{2},j}\) are the fluctuations across the boundary \(b_{i-\frac{1}{2},j}\) and where \(\mathcal{B}^{\pm}\Delta V_{i,j-\frac{1}{2}}\) are the fluctuations across the boundary \(b_{i,j-\frac{1}{2}}\). Inserting \(|V_{ij}|=r_i\Delta r \Delta \theta\) and \(l_{i-\frac{1}{2}}=r_{i-\frac{1}{2}} \Delta \theta\) into \eqref{eq:finiteVolumeScheme} and considering radially symmetric flow such that the fluctuations \(\mathcal{B}^+\Delta V_{i,j-\frac{1}{2}}\) and \(\mathcal{B}^-\Delta V_{i,j+\frac{1}{2}}\) vanish, we obtain the following numerical scheme for the simulation of the ASWME:
\begin{equation}\label{eq:finiteVolumeScheme2}
    V_{ij}^{n+1}=V_{ij}^n-\frac{\Delta t}{r_i\Delta r}\left( r_{i-\frac{1}{2}} \mathcal{A}^+\Delta V_{i-\frac{1}{2},j}+r_{i+\frac{1}{2}}\mathcal{A}^-\Delta V_{i+\frac{1}{2},j}\right).
\end{equation}
The fluctuations are computed using the PRICE-scheme \cite{PRICE} and the software framework also employed for the simulations in \cite{koellermeier_dissertation_2017}.

\subsubsection{Boundary conditions}

The physical domain together with a schematic overview of the used boundary conditions is illustarted in Figure \ref{fig:boundaryIllustration}. Boundary conditions need to be specified for the inner boundary (\(r=r_{min}\)) and for the outer boundary (\(r=r_{max}\)) while the angular direction is simply periodic. For the inner boundary, we impose a wall boundary condition, implemented as a reflective and impermeable boundary. We further assume that there is no slip along this boundary. 

Following the boundary condition treatment for solid walls in \cite{Steldermann_2023}, based on standard methods from \cite{LeVeque_2002}, the ghost cells are defined as

\begin{equation}
    \begin{pmatrix}
        h \\
        u_m \\
        \alpha_k \\
        v_m \\
        \gamma_k
    \end{pmatrix}^{\text{ghost}}
    =
    \begin{pmatrix}
        h^{\text{int., extrap.}} \\
        -u_m \\
        -\alpha_k \\
        0 \\
        0
    \end{pmatrix},
\end{equation}
where we use first order extrapolation to extrapolate the interior values for \(h\) to the ghost cell and where we changed the sign of the radial velocity (the normal velocity to the inner boundary). There is no slip along the wall, so the angular velocity \(v_m\) and the angular coefficients \(\gamma_k\) are set to zero. For the outer boundary, we impose an outflow boundary condition using zero-order extrapolation:

\begin{equation}
    \begin{pmatrix}
        h \\
        u_m \\
        \alpha_k \\
        v_m \\
        \gamma_k
    \end{pmatrix}^{\text{ghost}}
    =
    \begin{pmatrix}
        h \\
        u_m \\
        \alpha_k \\
        v_m \\
        \gamma_k
    \end{pmatrix}^{\text{int.}}.
\end{equation}

While the definition of more complex boundary conditions for the shallow water moment models is not fully solved in the literature, a rigorous analysis of the boundary conditions for moment models is beyond the scope of this paper and should be addressed in future work. In simulations, no visible problems with boundary effects were encountered.

\subsubsection{Reference solution}

To obtain a reference solution for the moment models, we follow \cite{SWME}, applied to a cylindrical framework. The numerical method for computing the reference solution is based on Equations \eqref{refsystemcylindrical1} - \eqref{refsystemcylindrical3}, but without the angular flux:

\vspace{-0.5cm}

\begin{align}\label{refsystemaxisymmetric1}
    &\partial_t h+\partial_r (hv_{r,m})+\frac{1}{r}hv_{r,m}=0,\\\label{refsystemaxisymmetric2}
    &\partial_t (hv_r)+\partial_r\left(hv_r^2+\frac{g}{2}e_zh^2\right)+\partial_{\zeta}\left(h v_r \omega-\frac{1}{\rho}\tau_{rz} \right)+\frac{h}{r}\left(v_r^2-v_{\theta}^2\right)=0,\\\label{refsystemaxisymmetric3}
    &\partial_t (hv_{\theta})+\partial_r (hv_rv_{\theta})+\partial_{\zeta}\left(h v_{\theta}\omega -\frac{1}{\rho}\tau_{\theta z} \right)+\frac{2h}{r}v_rv_{\theta}=0. 
\end{align}
The equations are split into a transport part that is solved explicitly and a source term part that is solved implicitly. For details, we refer to \cite{SWME}. The reference simulations have
been obtained using a grid with a spatial resolution of 200 cells in \(r\)-direction and 100 cells in \(\zeta\)-direction for a radial dam break in Section \ref{section:RadialDamBreakSimulation} and a grid with a spatial resolution of 400 cells in \(r\)-direction and 200 cells in \(\zeta\)-direction for a smooth test case in Section \ref{section:SmoothSimulation}, with a CFL number of 0.5. The approximations to the reference solution for both test cases showed sufficient convergence on the used grids. The reference solutions were obtained using a modified version of the software accompanying \cite{SWME}.

\subsection{Radial dam break}\label{section:RadialDamBreakSimulation}

\begin{figure}[h!]
    \centering
    \begin{subfigure}[b]{0.42\textwidth}
        \centering
        \includegraphics[width=\textwidth]{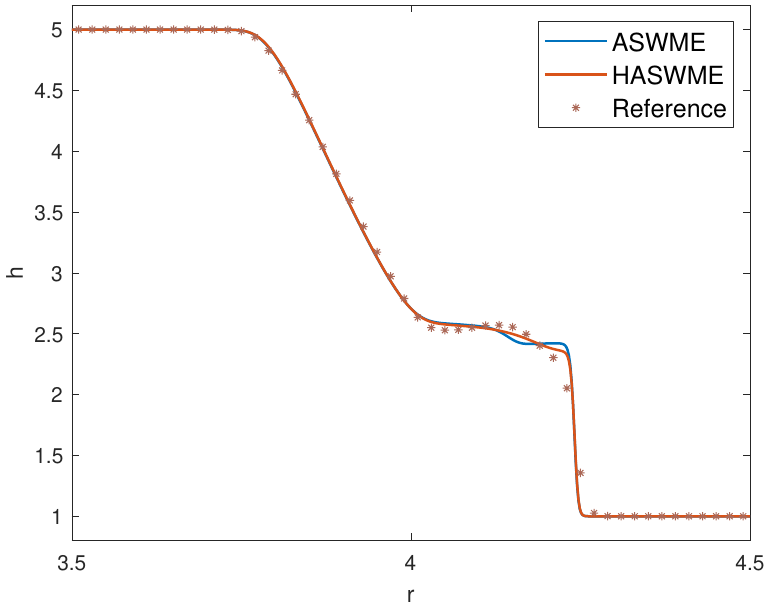}
        \caption{\(h\) at time \(t=0.1\)}
        \label{fig:radialDamBreak_h_Time0.1}
    \end{subfigure}
    \hfill
    \begin{subfigure}[b]{0.42\textwidth}
        \centering
        \includegraphics[width=\textwidth]{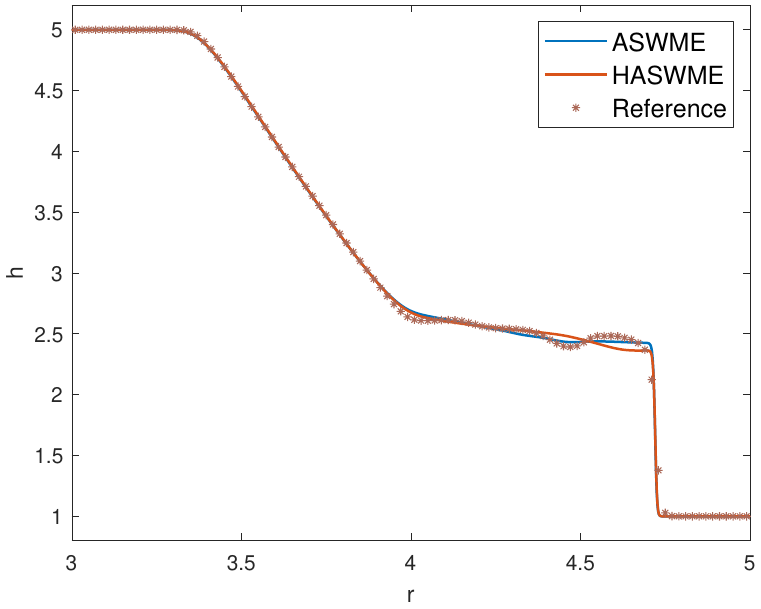}
        \caption{\(h\) at time \(t=0.3\)}
        \label{fig:radialDamBreak_h_Time0.3}
    \end{subfigure}
    \vskip\baselineskip
    \begin{subfigure}[b]{0.42\textwidth}
        \centering
        \includegraphics[width=\textwidth]{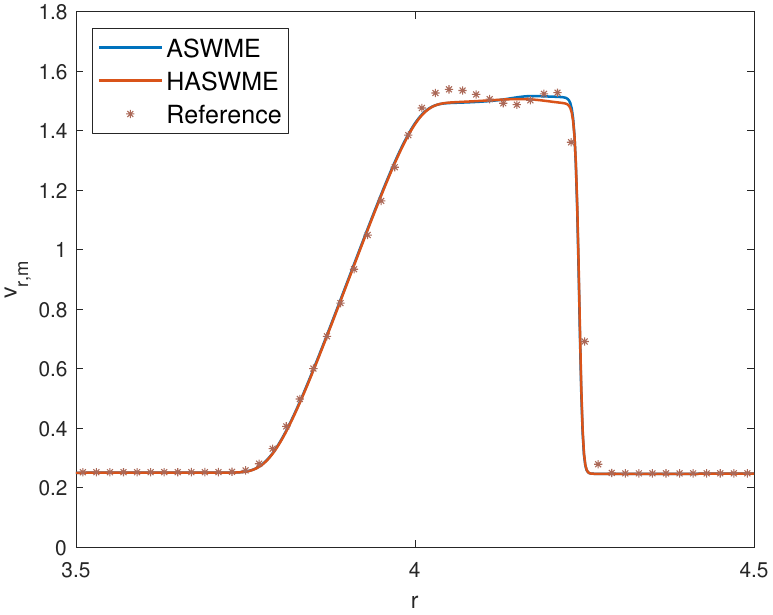}
        \caption{\(v_{r,m}\) at time \(t=0.1\)}
        \label{fig:radialDamBreak_vr_Time0.1}
    \end{subfigure}
    \hfill
    \begin{subfigure}[b]{0.42\textwidth}
        \centering
        \includegraphics[width=\textwidth]{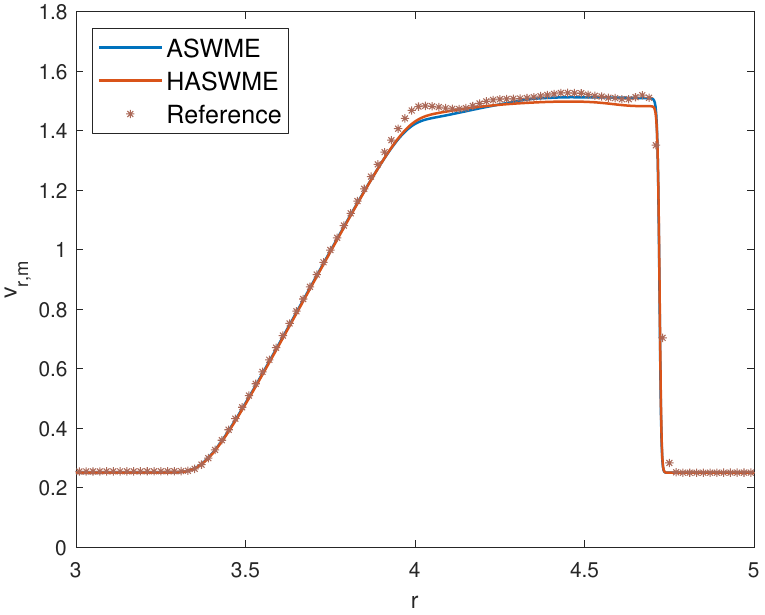}
        \caption{\(v_{r,m}\) at time \(t=0.3\)}
        \label{fig:radialDamBreak_vr_Time0.3}
    \end{subfigure}
    \vskip\baselineskip
    \begin{subfigure}[b]{0.42\textwidth}
        \centering
        \includegraphics[width=\textwidth]{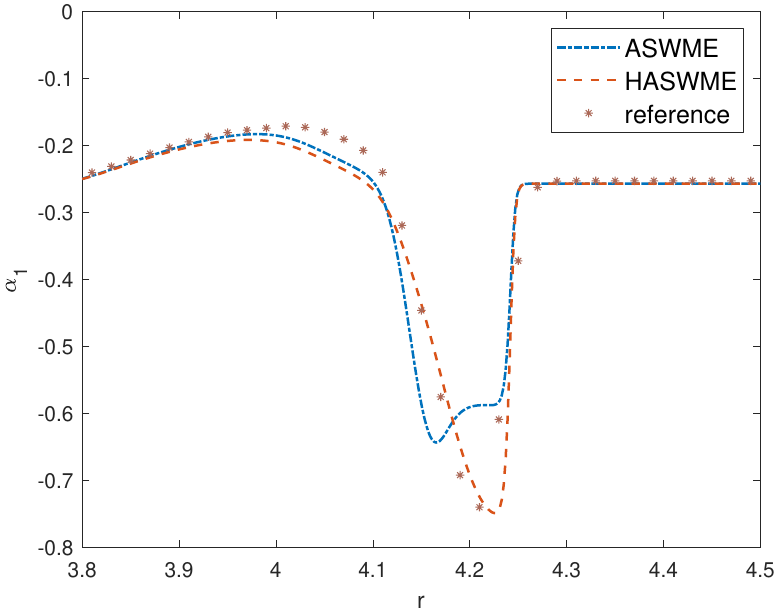}
        \caption{\(\alpha_1\) at time \(t=0.1\)}
        \label{fig:radialDamBreak_alpha1_Time0.1}
    \end{subfigure}
    \hfill
    \begin{subfigure}[b]{0.42\textwidth}
        \centering
        \includegraphics[width=\textwidth]{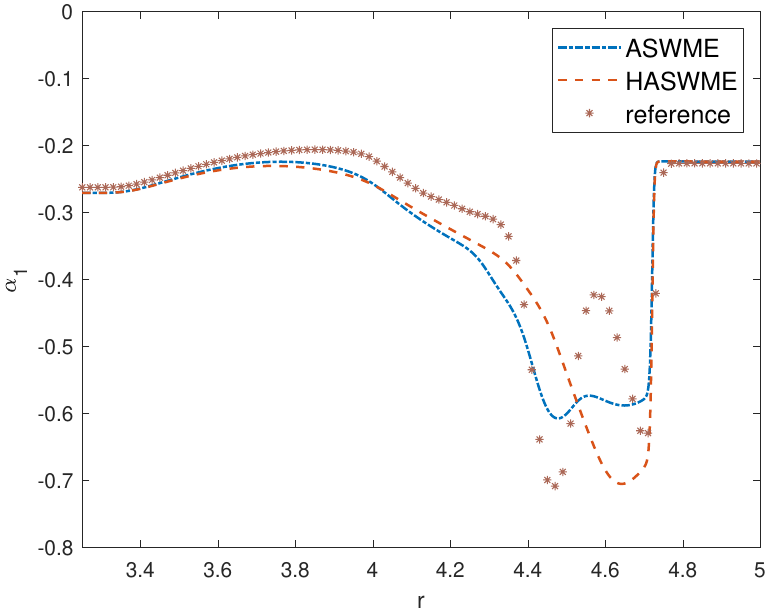}
        \caption{\(\alpha_1\) at time \(t=0.3\)}
        \label{fig:radialDamBreak_alpha1_Time0.3}
    \end{subfigure}
    \caption{Radial dam break test case for the third order ASWME and HASWME at times \(t=0.1\) and \(t=0.3\). HASWME gives comparable accuracy for \(h\) and \(v_{r,m}\) but improves the prediction of \(\alpha_1\) at \(t=0.1\).}
    \label{fig:RadialDamBreak}
\end{figure}

First, a test case with a discontinuous initial height profile will be considered. The numerical setup given in Table \ref{table:setup_RadialDamBreak} is very similar to the 1D test case with discontinuous initial data in \cite{HSWME}. The initial radial velocity is a cubic function of the water height. For this particular initial radial velocity, we expect hyperbolicity problems because of the observations in \cite{HSWME}. Note that in the situation where the angular velocity is zero everywhere the ASWME \eqref{systemFormCylindricalReduced} simply reduce to the 1D equations with additional geometric source terms, as the reference system \eqref{refsystemcylindrical1}-\eqref{refsystemcylindrical3} can then be reduced to simply \eqref{refsystemcylindrical1}-\eqref{refsystemcylindrical2}. To perform simulations of the pseudo 2D-model, we therefore set the initial angular velocity to be a linearly increasing function of the radial coordinate. The initial height function takes the form of a dam break in radial direction. The simulations are performed using the third order axisymmetric model denoted as ASWME and the third order hyperbolic axisymmetric model denoted as HASWME on an equidistant mesh with 2000 cells in radial direction. Simulation times are \(t=0.1\) and \(t=0.3\) and the results are plotted in Figure \ref{fig:RadialDamBreak}. The spatial domain \(r \in [2,6]\) is chosen such that significant influence of the geometric source term can be expected. The values of the variables are not plotted for the whole domain, but only in that part of the domain in which there are significant differences between the different approximations. 

\begin{table}[h!]
\begin{center}
\begin{tabular}{ | m{10em} || m{6.5cm} | } 
    \hline
    Kinematic viscosity & \(\nu=0.1\) \\[4pt] 
    \hline
    Slip length & \(\lambda=0.1\) \\[4pt] \hline 
    Spatial domain & \(r \in [2,6], \quad \theta \in [0,2\pi]\) \\[4pt]
    \hline
    End time & \(t_{end}\in \{0.1,0.3\}\) \\[4pt]
    \hline
    Inital height & \[
                    h(r,\theta,0) = 
                    \begin{cases}
                    5.0, &  r \leq 4.0 \\
                    1.0, &  r > 4.0
                    \end{cases}
                   \] \\[4pt] 
    \hline
    Initial radial velocity & \(v_r(r,\theta,\zeta,0)=0.25-2.5\zeta+7.5\zeta^2-5\zeta^3\) \\[4pt]
    \hline
    Initial angular velocity & \(v_{\theta}(r,\theta,\zeta,0)=0.1r\)\\[4pt]
    \hline
    Time integration & Forward Euler \\[4pt]
    \hline
    Spatial discretization & PRICE scheme \cite{PRICE} \\[4pt]
    \hline
    CFL number & 0.25 \\[4pt]
    \hline
\end{tabular}
\caption{Numerical setup for radial dam break scenario.}
\label{table:setup_RadialDamBreak}
\end{center}
\end{table}

%\enlargethispage{1\baselineskip}

\begin{figure}[h!]
    \centering
    \includegraphics[width=0.9\textwidth]{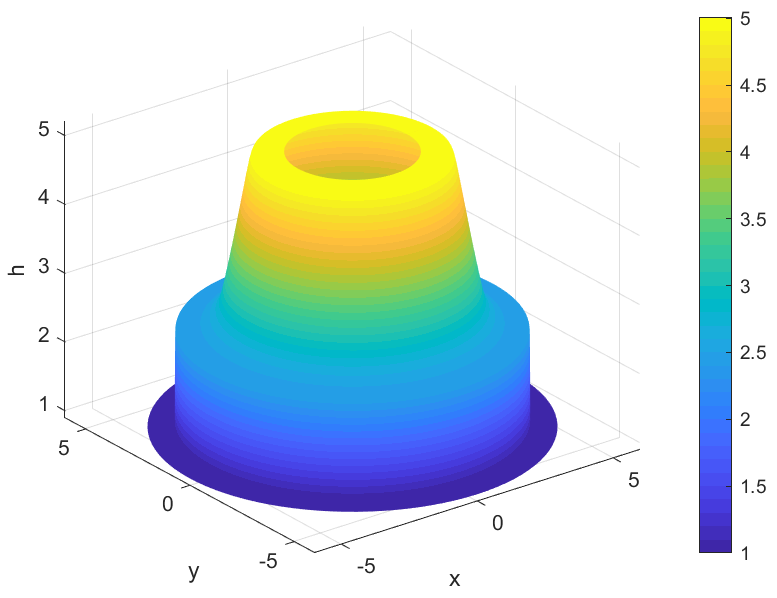}
    \caption{3D visualization of radial dam break simulation using the third order HASWME at time \(t=0.5\)}.
    \label{fig:3DRadialDam}
\end{figure}

First, a 3D plot of the simulation of the water height at time \(t=0.5\) is given in Figure \ref{fig:3DRadialDam}. This illustrates how the model is pseudo-two-dimensional: while the model equations are 1D and while the flow properties can be analyzed by considering slices in radial direction, the flow simulations can be represented in a 2D space. 

For this test case, Figure \ref{fig:radialDamBreak_h_Time0.1} and \ref{fig:radialDamBreak_h_Time0.3}, for the water height \(h\), and Figure \ref{fig:radialDamBreak_vr_Time0.1} and Figure \ref{fig:radialDamBreak_vr_Time0.3}, for the mean radial velocity \(v_{r,m}\), show that the ASWME and the HASWME yield similar and accurate results. For the approximation of the first coefficient \(\alpha_1\), however, the ASWME and the HASWME give qualitatively different results, as seen in Figure \ref{fig:radialDamBreak_alpha1_Time0.1} and \ref{fig:radialDamBreak_alpha1_Time0.3}. At time \(t=0.1\), the HASWME outperforms the ASWME, as the latter model displays an oscillation that is not present in the reference solution. This is in agreement with the numerical results obtained in \cite{HSWME}, in which an instability was observed in the numerical simulation of the 1D SWME. However, at time \(t=0.3\), which was not previously reported in the comparable test case in \cite{HSWME}, both the ASWME and the HASWME approximations are not accurate compared to the reference solution, so it is not straightforward to argue if one model is more accurate than the other. We note that both the reference solution and ASWME seem to predict an emerging instability in Figure \ref{fig:radialDamBreak_alpha1_Time0.3} while the HASWME model seems stable. This effect hints at increased stability of the hyperbolic model, similar as in \cite{HSWME} and should be further studied.

\subsection{Smooth initial height profile}\label{section:SmoothSimulation}
Next, we consider a test case with a smooth initial height profile. The numerical setup given in Table \ref{table:setup_Smooth} is closely related to the 1D smooth test case considered in \cite{HSWME}, but a slightly different initial height function and a different initial velocity profile are used. The initial height function takes the form of a sigmoid function. The initial radial velocity is zero everywhere, while the initial angular velocity is linearly increasing in the radial coordinate. We consider a spatial domain \(r\in [1,8]\) to simulate the waves caused by a small object falling into the water and to include significant effects of the geometric terms.

\begin{table}[h!]
\begin{center}
\begin{tabular}{ | m{10em} || m{4.5cm} | } 
    \hline
    Kinematic viscosity & \(\nu=1\) \\[4pt] 
    \hline
    Slip length & \(\lambda=0.1\) \\[4pt] \hline 
    Spatial domain & \(r\in [1,8], \quad \theta \in [0,2\pi]\) \\[4pt]
    \hline
    End time & \(t_{end}=1.0\) \\[5pt]
    \hline
    Initial height & 
                    \(h(r,\theta,0) = 3 - \frac{2}{1+e^{3(r-5)}}\) \\[4pt] 
    \hline
    Initial radial velocity & \(v_r(r,\theta,\zeta,0)=0\) \\[4pt]
    \hline
    Initial angular velocity & \(v_{\theta}(r,\theta,\zeta,0)=0.1 r\)\\[4pt]
    \hline
    Time integration & Forward Euler \\[4pt]
    \hline
    Spatial discretization & PRICE scheme \\[4pt]
    \hline
    CFL number & 0.25 \\[4pt]
    \hline
\end{tabular}
\caption{Numerical setup for smooth test case.}
\label{table:setup_Smooth}
\end{center}
\end{table}

\begin{figure}[ht!]
    \centering
    \begin{subfigure}[b]{0.42\textwidth}
        \centering
        \includegraphics[width=\textwidth]{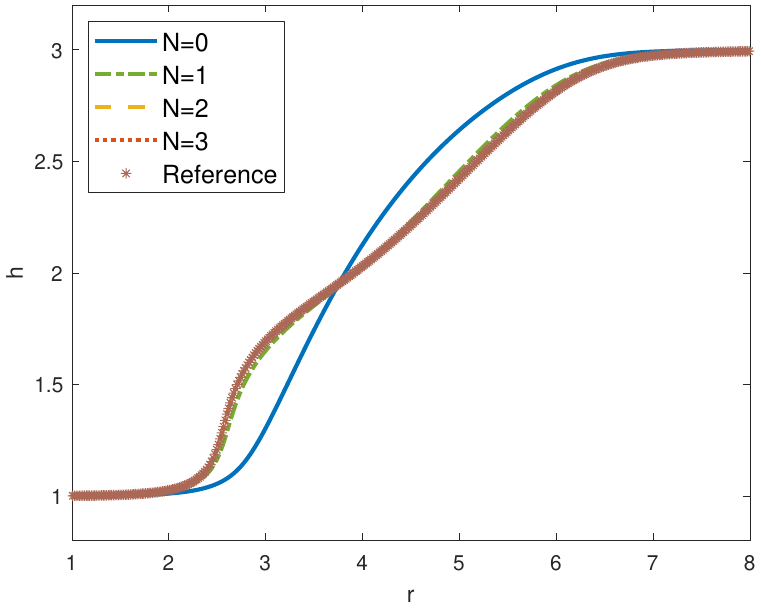}
        \caption{\(h\) at time \(t=1.0\)}
        \label{fig:smoothExpit_h_HASWME}
    \end{subfigure}
    \hfill
    \begin{subfigure}[b]{0.42\textwidth}
        \centering
        \includegraphics[width=\textwidth]{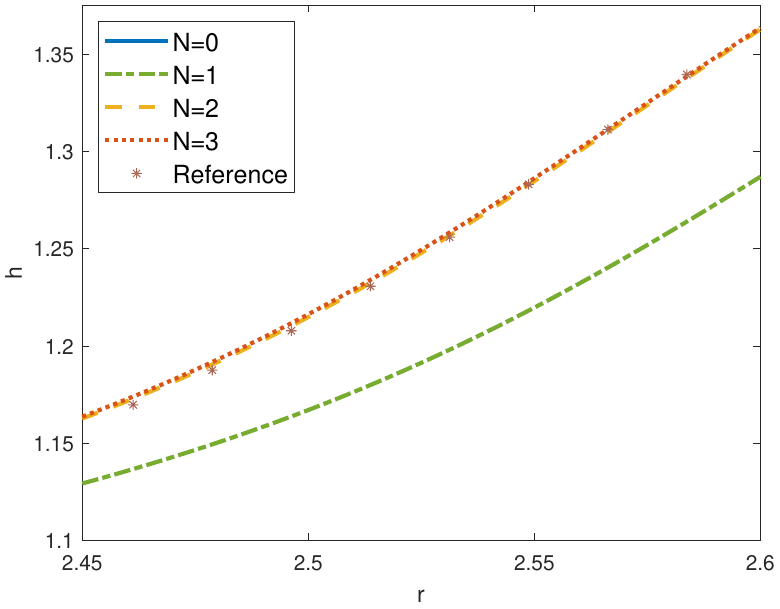}
        \caption{\(h\) at time \(t=1.0\)}
        \label{fig:smoothExpit_h_HASWME_ZoomedIn}
    \end{subfigure}
    \vskip\baselineskip
    \begin{subfigure}[b]{0.42\textwidth}
        \centering
        \includegraphics[width=\textwidth]{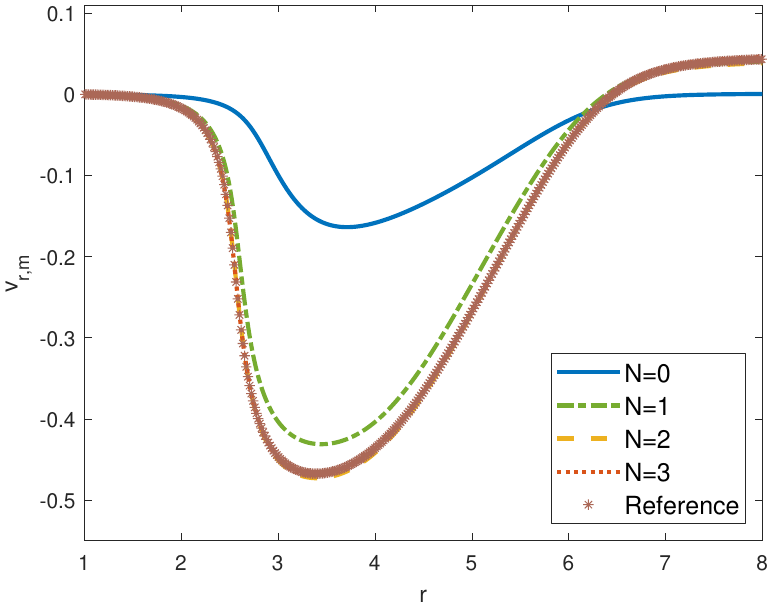}
        \caption{\(v_{r,m}\) at time \(t=1.0\)}
        \label{fig:smoothExpit_vr_HASWME}
    \end{subfigure}
    \hfill
    \begin{subfigure}[b]{0.42\textwidth}
        \centering
        \includegraphics[width=\textwidth]{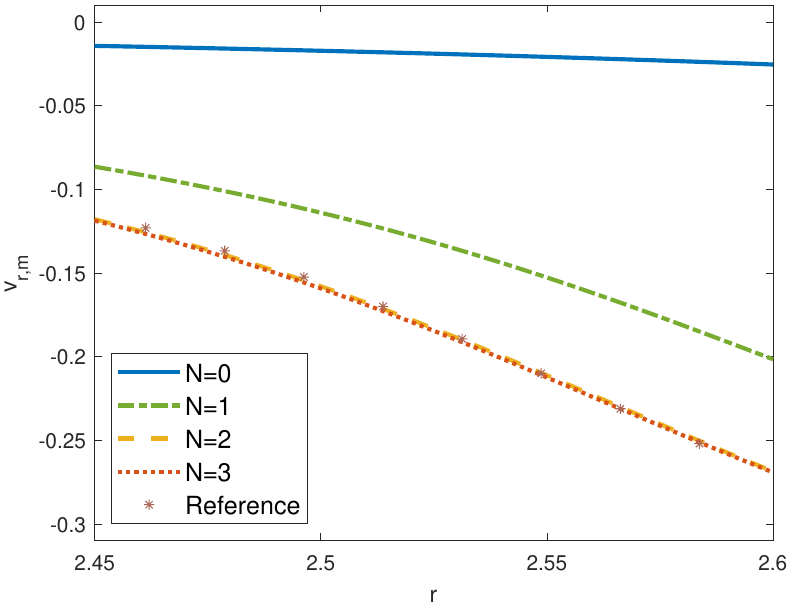}
        \caption{\(v_{r,m}\) at time \(t=1.0\)}
        \label{fig:smoothExpit_vr_HASWME_ZoomedIn}
    \end{subfigure}
    \vskip\baselineskip
    \begin{subfigure}[b]{0.42\textwidth}
        \centering
        \includegraphics[width=\textwidth]{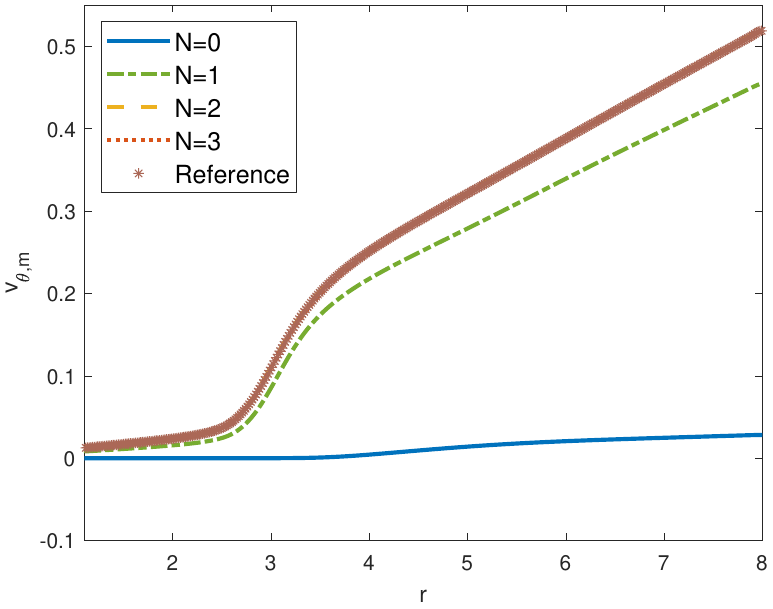}
        \caption{\(v_{\theta,m}\) at time \(t=1.0\)}
        \label{fig:smoothExpit_vtheta_HASWME}
    \end{subfigure}
    \hfill
    \begin{subfigure}[b]{0.42\textwidth}
        \centering
        \includegraphics[width=\textwidth]{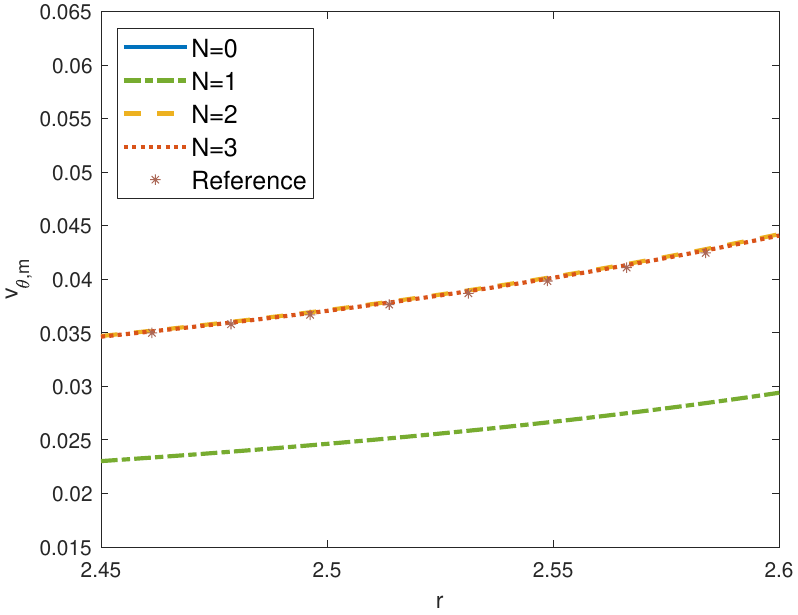}
        \caption{\(v_{\theta,m}\) at time \(t=1.0\)}
        \label{fig:smoothExpit_vtheta_HASWME_ZoomedIn}
    \end{subfigure}
    \caption{Smooth initial height profile test case for different orders \((N,N)\) HASWME at time \(t=1.0\) in the full domain (left column) and zoomed in on a segment of the radial domain (right column). The accuracy of the approximations increases with increasing order.}
    \label{fig:smoothExpit}
\end{figure}

Moment approximations are obtained using the HASWME \(N\)th order models, with \(N \in \{ 0,1,2,3 \}\), on an equidistant mesh with 4000 cells in radial direction. The results are shown in Figure \ref{fig:smoothExpit}. The simulation results for this smooth test case are very similar for the ASWME model, similar to Figure \ref{fig:smoothExpit}, as no issues of stability or hyperbolicity loss occur here. Hence the ASWME results are not shown here. The values of \(h\), \(v_{r,m}\) and \(v_{\theta,m}\) are plotted in the full domain (left column) and in a segment of the radial domain (right column) with increasing order $N$. Figure \ref{fig:smoothExpit_h_HASWME} shows that the approximations for the water height \(h\) are more accurate with increasing order, compared to the reference solution. This trend is also clearly visible for the mean radial velocity \(v_{r,m}\), see Figure \ref{fig:smoothExpit_vr_HASWME}, and for the mean angular velocity \(v_{\theta,m}\), displayed in Figure \ref{fig:smoothExpit_vtheta_HASWME}. The numerical values of the higher order models are so close to the reference solution that they can not be distinguished in Figures \ref{fig:smoothExpit_h_HASWME}, \ref{fig:smoothExpit_vr_HASWME} and \ref{fig:smoothExpit_vtheta_HASWME}. Therefore, Figures \ref{fig:smoothExpit_h_HASWME_ZoomedIn}, \ref{fig:smoothExpit_vr_HASWME_ZoomedIn} and \ref{fig:smoothExpit_vtheta_HASWME_ZoomedIn} zoom in on a segment of the radial domain. Clearly, the HASWME yield increasingly more accurate solutions with increasing number of moments, when compared to the reference solution. 

\begin{figure}[h!]
    \centering
    \begin{subfigure}[b]{0.47\textwidth}
        \centering
        \includegraphics[width=\textwidth]{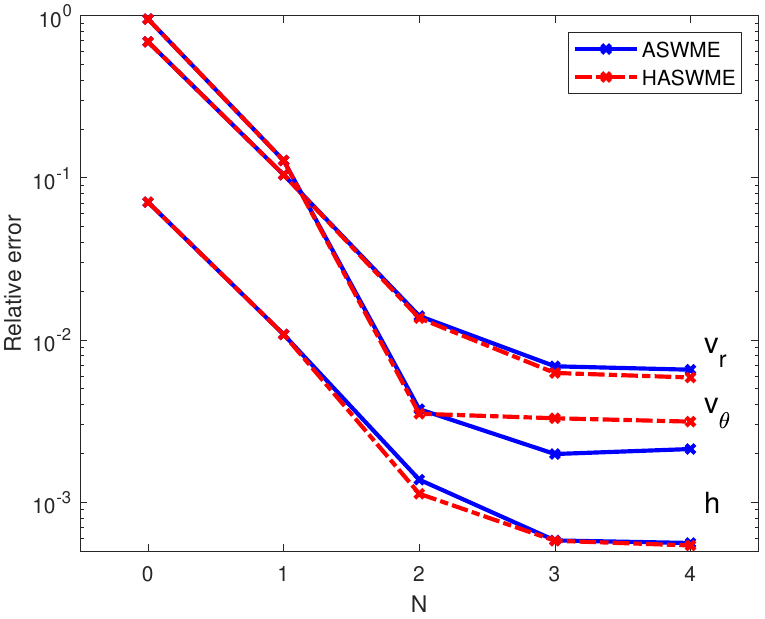}
        \caption{Fine simulation: 4000 grid cells for the moment simulations.}
        \label{fig:ErrorConvergenceFine}
    \end{subfigure}
    \hfill
    \begin{subfigure}[b]{0.47\textwidth}
        \centering
        \includegraphics[width=\textwidth]{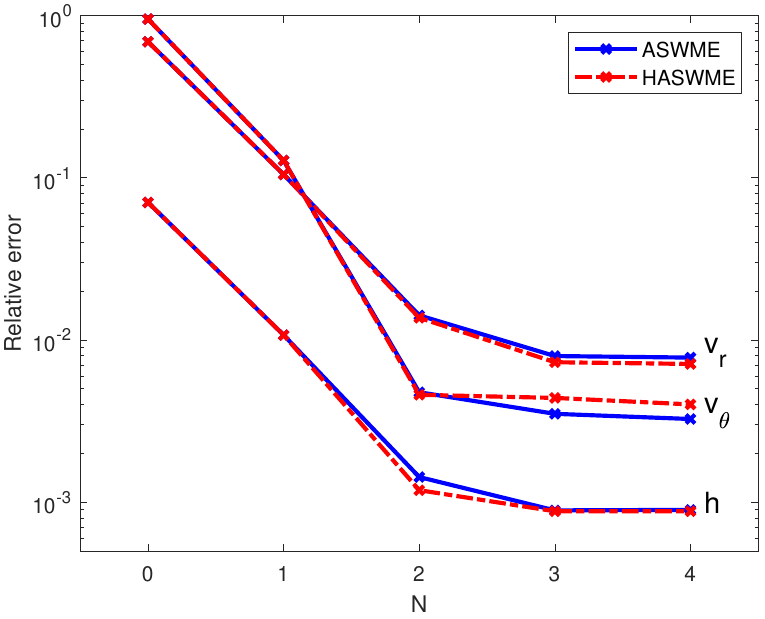}
        \caption{Coarse simulation: 2000 grid cells for the moment simulations.}
        \label{fig:ErrorConvergenceCoarse}
    \end{subfigure}
    \caption{Error convergence of smooth test case for ASWME and HASWME.}
    \label{fig:ErrorConvergence}
\end{figure}

The error convergence is shown in Figure \ref{fig:ErrorConvergence}, which shows the relative error of the different models for the water height \(h\), the mean radial velocity \(v_{r,m}\) and the mean angular velocity \(v_{\theta,m}\) for both the ASWME and the HASWME on different grids for the moment model simulations and using a grid with 400 grid cells in radial direction \(r\) and 200 grid cells in vertical direction \(\zeta\) in the reference simulation. Figure \ref{fig:ErrorConvergenceFine} shows the numerical results obtained using a grid with 4000 grid cells in radial direction \(r\) in the moment model simulations. For all three variables, we observe a reduction of the error when the order increases from \(N=0\) to \(N=3\). Note that ASWME and HASWME are identical for \(N=0\) and \(N=1\). In particular, the error in the ASWE model ($N=0$) is considerably larger than the error in both moment models ASWME and HASWME for $N>0$. This is in agreement with the 1D results in \cite{HSWME}. Note that when the order is increased from \(N=3\) to \(N=4\), there is no reduction in the error anymore for most of the models. However, the error is already small in the third order model, especially for \(h\) and \(v_{\theta,m}\), and the simulations are performed on finite grids leading to numerical errors that can dominate the error when the modelling error is small. To illustrate that, Figure \ref{fig:ErrorConvergenceCoarse} shows the results obtained using a coarser grid for the moment model simulations while keeping the CFL number constant. The moment model simulations use a grid with 2000 cells in radial direction \(r\). There is a smaller reduction in the relative error for the water height \(h\), the mean radial velocity \(v_{r,m}\) and the mean angular velocity \(v_{\theta,m}\) compared to the finer grid moment model simulations shown in Figure \ref{fig:ErrorConvergenceFine}. The same stagnation of the relative error when the order is increased from \(N=3\) to \(N=4\) can be observed for both grids. For these orders \(N=3\) and \(N=4\), the relative error on the finer grid, shown in Figure \ref{fig:ErrorConvergenceFine}, is visibly smaller than the relative error on the coarser grid, shown in Figure \ref{fig:ErrorConvergenceCoarse}. This implies that the relative errors on the two grids shown in Figure \ref{fig:ErrorConvergence} are dominated by the numerical error and that the modeling error is already relatively small. 

We can conclude that for this smooth test case, the moment models yield approximations that are increasingly accurate when the order of the model increases. This is the case for both the non-hyperbolic and the hyperbolic models. Furthermore, the HASWME yield similar or sometimes better accuracy compared to ASWME, even though it has used the hyperbolic regularization. Finally, the differences between ASWME and HASWME are very small compared to the overall error and the \(N=0\) model.

\clearpage

%% file: axisymmetric/Conclusion.tex
\section{Conclusion}
In this paper we derived, analyzed and numerically solved the first hyperbolic moment models for radially symmetric shallow flow in cylindrical coordinate systems. The paper started with a brief review of the derivation of the SWME. The reference system was then transformed to cylindrical coordinates and Axisymmetric Shallow Water Moment Equations (ASWME) were formulated for radially symmetric flow. We showed the breakdown of hyperbolicity of the ASWME by means of hyperbolicity plots. This lead to the derivation of the first hyperbolic axisymmetric model called Hyperbolic Axisymmetric Shallow Water Moment Equations (HASWME) for which we explicitly constructed the system matrix for an arbitrary number of moments and analyzed the eigenvalues of the system. A finite volume scheme tailored to the cylindrical geometry was derived for the ASWME. In the numerical simulation of a test case with a radial dam break, the hyperbolic model resulted in a more accurate approximation compared to a reference solution. The numerical simulation of a test case with smooth initial data showed that the approximation error converges when the order of the model is increased. We conclude that the new HASWME is a successful hyperbolic regularization of the ASWME with high accuracy and convergence towards the reference solution, while exhibiting better stability properties due to their hyperbolicity. Although a lot of the structure of the equations is inherited from the 1D equations, The models derived in this paper allow for the first simulation of 2D axisymmetric follow using hyperbolic moment equations. The models can therefore be seen as an intermediate stage between the 1D and 2D models.

Ongoing work focuses on the extension of the axisymmetric models to moment models that are fully 2D and the numerical simulation of these models with this paper paving the way for an extension of the analysis, especially for hyperbolicity. Suggestions for future work are the equilibrium stability analysis of the ASWME and the HASWME and the investigation of more advanced test cases.

\section*{Data Availability Statement}
The datasets generated and analysed during this study are available from the corresponding author on reasonable request.

\section*{Acknowledgements}
The authors would like to acknowledge the financial support of the CogniGron research center and the Ubbo Emmius Funds (University of Groningen). 

%% file: axisymmetric/Appendix.tex
\appendix

\section{Analytical construction of the hyperbolic system matrix}

\label{ch:hyperbolicityproof(N,N)}

\begin{reptheorem}{NNOrderTheoremMatrixForm}
    The HASWME system matrix \(A_{HA}^{(N,N)} \in \mathbb{R}^{(2N+3)\times (2N+3)}\) is given by

    \begin{equation*}
        A_{HA}^{(N,N)}=
        \begin{pmatrix}
            \boldsymbol{A}^{(N,N)} & \boldsymbol{0}^{(N,N)} \\
            \boldsymbol{B}^{(N,N)} & \boldsymbol{C}^{(N,N)}
        \end{pmatrix},
    \end{equation*}
    with \(\boldsymbol{A}^{(N,N)} \in \mathbb{R}^{(N+2)\times (N+2)}\) given by \eqref{ANN}, \(\boldsymbol{B}^{(N,N)}\in \mathbb{R}^{(N+1)\times (N+2)}\) given by \eqref{BNN} and \(\boldsymbol{C}^{(N,N)} \in \mathbb{R}^{(N+1)\times (N+1)}\) given by \eqref{CNN} and where all other entries are zero. \(\boldsymbol{0}^{(N,N)} \in \mathbb{R}^{(N+2)\times (N+1)}\) is a zero matrix.
\end{reptheorem}

\begin{proof}
The matrix \(A_{HA}^{(N,N)}\) is obtained by computing \(A_A^{(N,N)}\), the full system matrix, and then setting the higher order moments to zero \cite{HSWME}. We derive the rows of the system matrix separately:

\textbf{1. Mass and radial momentum balance equations.}
The first two rows of the system matrix correspond to the mass balance and the radial momentum balance equations. It can be easily seen that the first two rows of the coefficient matrix are given by
\begin{align*}
    \begin{pmatrix}
        & 1 & & & & & & & \\[6pt]
        gh-v_{r,m}^2-\frac{1}{3}\alpha_1^2+\sum_{i=1}^{N}\alpha_i^2 & 2v_{r,m} & \frac{2}{3} \alpha_1 & \frac{2}{5}\alpha_2 & \cdots & \frac{2}{2N+1}\alpha_N & 0 & \cdots & 0
    \end{pmatrix}.
\end{align*}
After setting the higher order moments to zero, we obtain
\begin{align*}
    \begin{pmatrix}
        & 1 & & & & & \\[6pt]
        gh-v_{r,m}^2-\frac{1}{3}\alpha_1^2 & 2v_{r,m} & \frac{2}{3} \alpha_1 &  &  &  & 
    \end{pmatrix}
\end{align*}
as the first to rows of \(A_{HA}^{(N,N)}\).

\textbf{2. The equations for \(\boldsymbol{\alpha_i,i=1,\ldots,N}\)}.
The equation for the \(i\)th moment reads
\begin{equation*}
    \frac{\partial(h\alpha_i)}{\partial t}+\frac{\partial F_i}{\partial r}=Q_i+P_i.
\end{equation*}
The system matrix \(A_H^{(N,N)}\) is composed of a conservative and a non-conservative part.

\textbf{2.1 The conservative part \(\boldsymbol{\frac{\partial F_i}{\partial r}}\).} Recall that
\begin{align*}
    F_i&=h\left( 2v_{r,m}\alpha_i+\sum_{j,k=1}^{N}A_{ijk}\alpha_j\alpha_k \right)=h2v_{r,m}\alpha_i+h\sum_{j,k=1}^{N}(2i+1)\left(\int_0^1\phi_i \phi_j\phi_k d\zeta\right) \alpha_j\alpha_k.
\end{align*}
First note that this term does not depend on the angular moments, so all the partial derivative with respect to these moments will be zero. Consider the first term, \(h2v_{r,m}\alpha_i\). With \( \boldsymbol{\alpha}=(\alpha_1,\cdots,\alpha_{N})\), it can be easily seen (after setting the higher coefficients to zero) that this term leads to
\begin{equation*}
    \frac{\partial (2hv_{r,m}\boldsymbol{\alpha})}{\partial V}\quad\myeq\quad\begin{pmatrix}
        -2v_{r,m}\alpha_1 & 2\alpha_1 & 2v_{r,m} & & & & \\[6pt]
        & & & \ddots & & & \\[6pt]
        & & & & 2v_{r,m} &  &
    \end{pmatrix},   
\end{equation*}

where \(\quad \myeq \quad\) is setting the higher order moments to zero.

Now, we look at the second term, with the triple Legendre integral. We call this term \(\Tilde{F_i}:=h\sum_{j,k=1}^{N}A_{ijk}\alpha_j\alpha_k\). In \cite{HSWME}, it is shown that if \(j=1\), the only cases that need to be considered are \(k=i-1\) (if \(i\geq 2\)) and \(k=i+1\). Analogously, if \(k=1\), the only cases we need to consider are \(j=i-1\) (if \(i\geq 2\)) and \(j=i+1\). According to this observations, we %only need to look at a few combinations for the index triplet \((i,j,k)\). We 
can distinguish four different situations.

\textbf{2.1.1 Case 1: \(i=1\)}.
If \(i=1\), there are two cases for the index triplet \((i,j,k)\), with \(j=1\) or \(k=1\), which lead to a non-zero term: $(1,1,2) \text{ and } (1,2,1)$.

So we can pull these cases outside of the sum and we get:
\begin{equation*}
    \Tilde{F_1}=\sum_{j,k=1}^{N}A_{1jk}\alpha_j\alpha_k=2A_{112}h\alpha_1\alpha_2+h\sum_{j,k=2}^{N}A_{1jk}\alpha_j\alpha_k.
\end{equation*}
The factor \(2\) appears because both index triplets lead to the same term. 

The derivative of the second term with respect to $h$ reads
\begin{equation*}
    \frac{\partial \left(h\sum_{j,k=2}^{N}A_{1jk}\alpha_j\alpha_k\right)}{\partial h}=-\sum_{j,k=2}^{N}A_{1jk}\alpha_j\alpha_k.
\end{equation*}
After setting the higher order moments to zero, this term reduces to zero. The same can be observed when taking the derivative with respect to \(hv_{r,m}\), \(h\alpha_l\) (with \(l=1,\cdots,N\)) and \(hv_{\theta,m}\). 

The derivatives of the first term read
\begin{equation*}
    \frac{\partial(2A_{112}h\alpha_1\alpha_2)}{\partial(h\alpha_l)}=2A_{112}\alpha_1 \delta_{2,l}+2A_{112}\alpha_2 \delta_{1,l}.
\end{equation*}

When we set higher order moments to zero, the second term vanishes. We call the resulting entry \(c_1^F\). This entry corresponds to the derivative with respect to \(h\alpha_2\) and will be discussed later.
All other derivatives are zero, except the derivative with respect to \(h\), but we can easily see that this derivative reduces to zero when we set higher order moments to zero.

\textbf{2.1.2 Case 2: \(i=2\)}. 
There are three cases for the index triplet \((i,j,k)\), with \(j=1\) or \(k=1\), which lead to a non-zero term: $(2,1,1), (2,1,3), \text{ and } (2,3,1)$.

We can pull these terms outside of the sum:
\begin{equation*}
    \Tilde{F_2}=\sum_{j,k=1}^{N}A_{2jk}\alpha_j\alpha_k=A_{211}h\alpha_1^2+2A_{213}h\alpha_1\alpha_3+h\sum_{j,k=2}^{N}A_{2jk}\alpha_j\alpha_k.
\end{equation*}

Again, the derivatives of the last term all reduce to zero when setting the higher order moments to zero. Denoting the first and the second term by \(\underline{\Tilde{F}_2}:=A_{211}h\alpha_1^2+2A_{213}h\alpha_1\alpha_3\), we have:
\begin{align*}
    \frac{\partial\underline{\Tilde{F}_2}}{\partial (h\alpha_l)}=2A_{211}\alpha_1\delta_{1,l}+2A_{213}\alpha_1\delta_{3,l}+2A_{213}\alpha_3\delta_{1,l}.
\end{align*}
Setting higher order moments to zero, this reduces to 
\begin{equation*}
    2A_{211}\alpha_1\delta_{1,l}+2A_{213}\alpha_1\delta_{3,l}.
\end{equation*}
The first term leads to an entry \(a_2^F\). The second term leads to an entry \(c_2^F\).
In the same way, one can easily see that
\begin{equation*}
    \frac{\partial\underline{\Tilde{F}_2}}{\partial h}=-A_{211}\alpha_1^2-\alpha_1\alpha_3,\quad\frac{\partial\underline{\Tilde{F}_2}}{\partial (hv_{r,m})}=0,\quad\frac{\partial\underline{\Tilde{F}_2}}{\partial (hv_{\theta,m})}=0.
\end{equation*}
When higher order moment are set to zero, this only results in a term \(-\frac{2}{3}\alpha_1^2\) in the first column.
\vspace{3mm}

\textbf{2.1.3 Case 3: \(2<i\leq N-1\)}.
In this case, there are four possibilities for the index triplet \((i,j,k)\): $(i,1,i-1), (i,1,i+1), (i,i-1,1), \text{ and } (i,i+1,1)$.

We can pull this terms outside of the sum:
\begin{equation*}
    \Tilde{F_i}=\sum_{j,k=1}^{N}A_{ijk}\alpha_j\alpha_k=2A_{i,1,i-1}h\alpha_1\alpha_{i-1}+2A_{i,1,i+1}h\alpha_1\alpha_{i+1}+h\sum_{j,k=2}^{N}A_{ijk}\alpha_j\alpha_k.
\end{equation*}

The derivatives of the last term all reduce to zero when setting the higher order moments to zero. Denoting the remainder by \(\underline{\Tilde{F}_i}:=2A_{i,1,i-1}h\alpha_1\alpha_{i-1}+2A_{i,1,i+1}h\alpha_1\alpha_{i+1}\) yields
\begin{equation*}
    \frac{\partial\underline{\Tilde{F}_i}}{\partial (h\alpha_l)}=2A_{i,1,i-1}\alpha_1\delta_{i-1,l}+2A_{i,1,i-1}\alpha_{i-1}\delta_{1,l} +2A_{i,1,i+1}\alpha_1\delta_{i+1,l}+2A_{i,1,i+1}\alpha_{i+1}\delta_{1,l}.
\end{equation*}
%\begin{multline*}
%    \frac{\partial\underline{\Tilde{F}_i}}{\partial (h\alpha_l)}=2A_{i,1,i-1}\alpha_1\delta_{i-1,l}+2A_{i,1,i-1}\alpha_{i-1}\delta_{1,l}\\+2A_{i,1,i+1}\alpha_1\delta_{i+1,l}+2A_{i,1,i+1}\alpha_{i+1}\delta_{1,l}.
%\end{multline*}
Setting higher order moments to zero, this reduces to 
\begin{equation*}
    2A_{i,1,i-1}\alpha_1\delta_{i-1,l}+2A_{i,1,i+1}\alpha_1\delta_{i+1,l}.
\end{equation*}
The first term leads to an entry \(a_i^F\). The second term leads to an entry \(c_i^F\).
For the other derivatives, it is again easy to see that
\begin{equation*}
    \frac{\partial\underline{\Tilde{F}_i}}{\partial h}=-2A_{i,1,i-1}\alpha_1\alpha_{i-1}-\alpha_1\alpha_{i+1},\quad\frac{\partial\underline{\Tilde{F}_i}}{\partial (hv_{r,m})}=0,\quad\frac{\partial\underline{\Tilde{F}_i}}{\partial (hv_{\theta,m})}=0.\\
\end{equation*}
When higher order moments are set to zero, the derivative with respect to \(h\) reduces to zero.

\textbf{2.1.4 Case 4: \(i=N\)}.
In this case, there are two possibilities for the index triplet \((i,j,k)\): $(N,1,N-1) \text{ and } (N,N-1,1)$.

Analogously to the previous cases, it can be easily shown that this equation only leads to an entry \(a_{N}^F\) in the \((N-1)\)th column of the \((N+2)\)th row, the row corresponding to the equation for the \(N\)th moment.

Using orthogonality and recursion formulas (see also \cite{HSWME}), we have:
\begin{align*}
    &a_i^F=\frac{2i}{2i-1}\alpha_1, \qquad i=2,\ldots,N,\\[8pt]
    &c_i^F=\frac{2i+2}{2i+3}\alpha_1, \qquad i=1,\ldots,N-1.
\end{align*}
The entry \(a_i^F\) corresponds to a derivative with respect to \(h\alpha_{i-1}\) and is located on the first lower diagonal. The entry \(c_i^F\) corresponds to a derivative with respect to \(h\alpha_{i+1}\) and is located on the first upper diagonal leading to the modified flux derivative \(\frac{\partial F}{\partial V}_{mod}\):
\begin{equation}
    \frac{\partial F}{\partial V}_{mod}=
    \begin{pmatrix}
        & 1 & & & & &  & & \\[6pt]
        gh-v_{r,m}^2-\frac{1}{3}\alpha_1^2 & 2v_{r,m} & \frac{2}{3}\alpha_1 & & & & & \\[6pt]
        -2 v_{r,m}\alpha_1 & 2\alpha_1 & 2v_{r,m} & c_1^F & & & & & \\[6pt]
        -\frac{2}{3}\alpha_1^2 & & a_2^F & 2v_{r,m} & \ddots & & & & \\[6pt]
        & & & \ddots & \ddots & c_{N-1}^F & & \\[6pt]
        & & & & a_{N}^F & 2v_{r,m} & & \\[6pt]
        -v_{r,m}v_{\theta,m} & v_{\theta,m} & & & & & v_{r,m} & 
    \end{pmatrix},
\end{equation}
where the entries corresponding to partial derivatives with respect to the coefficients \(h\gamma_l\) are zero.

\textbf{2.2 The non-conservative part \(\boldsymbol{Q_i}\)}.
Recall the non-conservative part in the equation for \(h\alpha_i\):
\begin{equation}\label{proof_Nr_nc}
    Q_i=v_{r,m}\frac{\partial(h\alpha_i)}{\partial r}-\sum_{j,k=1}^{N}B_{ijk}\alpha_k\frac{\partial(h\alpha_j)}{\partial r}.
\end{equation}
Again, this term does not depend on the angular basis coefficients \(h\gamma_l\), so the partial derivatives of \(Q_i\) with respect to these terms will be zero. Clearly, the first term leads to the term \(v_{r,m}\) on the diagonal of \(Q\). The second term is simplified when we set higher order moments to zero:
\begin{equation*}
    -\sum_{j=1}^{N}(2i+1)\left(\int_0^1\phi_i'\left( \int_0^1\phi_jd\hat{\zeta}\right)\phi_1d\zeta\right)\alpha_1\frac{\partial(h\alpha_j)}{\partial r}.
\end{equation*}
It can be proved that \(B_{ij1}=0\) for \(|i-j|\neq 1\), except for the first two rows \cite{HSWME}. \(Q_i\) does not depend on any derivative with respect to \(h,hv_{r,m}\) and \(hv_{\theta,m}\), so the entries of the corresponding columns will be zero. The entries of the first two rows and the last row are also all zero; the mass and momentum balances do not contain non-conservative terms. By considering three different cases (\(i=1\), \(1<i\leq N-1\) and \(i=N\)) in an analogous way as for the conservative terms, we can see that the second term in \eqref{proof_Nr_nc} only leads to non-zero entries \(a_i^Q\) on the first lower diagonal and non-zero entries \(c_i^Q\) on the first upper diagonal:
\begin{align*}
    &a_i^Q=\frac{i+1}{2i-1}\alpha_1 \qquad i=2,\ldots,N,\\[8pt]
    &c_i^Q=\frac{i}{2i+3}\alpha_1 \qquad i=1,\ldots,N-1.
\end{align*}

This results in the following non-conservative part for the first \(N+2\) equations: 

\begin{equation}
Q_{N+2,mod}=\begin{pmatrix} Q_{(N+2)\times (N+3)} & 0_{(N+2)\times N}\end{pmatrix}\in \mathbb{R}^{(N+2)\times (2N+3)},
\end{equation}
with \(0_{(N+2)\times N} \in \mathbb{R}^{(N+2)\times N}\) a zero matrix and with
\begin{equation}\label{nonconsflux}
    Q_{(N+2)\times (N+3)}=
    \begin{pmatrix}
        & & & & & & & \\[6pt]
        & & & & & & & \\[6pt]
        & & & v_{r,m} & c_1^Q & & & \\[6pt]
        & & & a_2^Q & v_{r,m} & \ddots & & \\[6pt]
        & & & & \ddots & \ddots & c_{N-1}^Q & \\[6pt]
        & & & & & a_{N}^Q & v_{r,m} & \\[6pt]
        & & & & & & & 
    \end{pmatrix} \in \mathbb{R}^{(N+2)\times (N+3)}.
\end{equation}

Putting the conservative contributions  and the non-conservative contributions together, we obtain

\begin{equation*}
    \begin{pmatrix}
        & 1 & 0 & \cdots & & & & &  0  \\[6pt]
        gh-v_{r,m}^2-\frac{1}{3}\alpha_1^2 & 2v_{r,m} & \frac{2}{3}\alpha_1 & 0 & \cdots & & & & 0\\[6pt]
        -2 v_{r,m}\alpha_1 & 2\alpha_1 & v_{r,m} & \frac{3}{5}\alpha_1 & 0 & \cdots & & & 0\\[6pt]
        -\frac{2}{3}\alpha_1^2 & & \frac{1}{3}\alpha_1 & v_{r,m} & \ddots & 0 & \cdots &  & 0\\[6pt]
        & & & \ddots & \ddots & \frac{N+1}{2 N+1}\alpha_1 & 0 & \cdots & 0\\[6pt]
        & & & & \frac{N-1}{2N-1}\alpha_1 & v_{r,m} & 0 & \cdots & 0
    \end{pmatrix}
\end{equation*}
as the first \(2+N\) rows of the matrix \(A_{HA}^{(N,N)}\).

\textbf{3. Angular momentum balance equation.}
The \((N+3)\)th row of the coefficient matrix corresponds to the angular momentum balance equation. Recall that the flux in this equation reads
\begin{equation*}
    F_{N+3}:=hv_{r,m}v_{\theta,m}+h\sum_{j=1}^{N}\frac{\alpha_j\gamma_j}{2j+1}.
\end{equation*}
It can be easily verified that
\begin{align*}
    &\frac{\partial F_{N+3}}{\partial h}=-v_{r,m}v_{\theta,m}-\sum_{j=1}^{N}\frac{\alpha_j\gamma_j}{2j+1}\quad\myeq\quad -v_{r,m}v_{\theta,m}-\frac{\alpha_1\gamma_1}{3},\\[6pt]&\frac{\partial F_{N+3}}{\partial (hv_{r,m})}=v_{\theta,m}\quad\myeq\quad v_{\theta,m},\quad\frac{\partial F_{N+3}}{\partial (h\alpha_l)}=\frac{\gamma_l}{2l+1}\quad\myeq\quad \frac{\gamma_l}{2l+1}\delta_{1,l},\\[6pt]&\frac{\partial F_{N+3}}{\partial (hv_{\theta,m})}=v_{r,m}\quad\myeq\quad v_{r,m},\quad
    \frac{\partial F_{N+3}}{\partial (h\gamma_l)}=\frac{\alpha_l}{2l+1}\quad\myeq\quad \frac{\alpha_l}{2l+1}\delta_{1,l}.
\end{align*}

\textbf{4. The equations for \(\boldsymbol{\gamma_i,i=1,\ldots,N}\)}.
The equation for the \(i\)th moment is given by
\begin{equation*}
    \frac{\partial(h\gamma_i)}{\partial t}+\frac{\partial F_i}{\partial r}=Q_i+P_i.
\end{equation*}

The system matrix \(A_{HA}^{(N,N)}\) is composed of a conservative and a non-conservative part.

\textbf{4.1 The conservative part \(\boldsymbol{\frac{\partial F_i}{\partial r}}\)}.
Recall that
\begin{align}
    F_i&=h\left( v_{r,m}\gamma_i+v_{\theta,m}\alpha_i+\sum_{j,k=1}^{N}A_{ijk}\alpha_j\gamma_k \right)\\
    &=hv_{r,m}\gamma_i+hv_{\theta,m}\alpha_i+h\sum_{j,k=1}^{N}(2i+1)\left(\int_0^1\phi_i \phi_j\phi_k d\zeta\right) \alpha_j\gamma_k. \label{consflux}
\end{align}
Consider the first term in the right-hand side in Equation \eqref{consflux},  \(hv_{r,m}\gamma_i\). Clearly, this term does not depend on \(\alpha_l\) (\(l=1,\ldots,N\)) and on \(v_{\theta,m}\), so all the derivatives with respect to \(h\alpha_l\) and \(v_{\theta,m}\) will be zero. Furthermore, we have:
\begin{equation*}
    \frac{\partial (hv_{r,m}\gamma_i)}{\partial h}=-v_{r,m}\gamma_i, \quad\frac{\partial (hv_{r,m}\gamma_i)}{\partial (hv_{r,m})}=\gamma_i, \quad\frac{\partial (hv_{r,m}\gamma_i)}{\partial (h\gamma_l)}=v_{r,m}\delta_{i,l}.
\end{equation*}
Then we consider the second term in the right-hand side in Equation \eqref{consflux}, \(hv_{\theta,m}\alpha_i\). This term does not depend on \(\gamma_l\) (\(l=1,\ldots,N\)) and on \(hv_{\theta,m}\), so all the derivatives with respect to \(h\gamma_l\) and with respect to \(hv_{\theta,m}\) will be zero. Furthermore, we have:
\begin{equation*}
    \frac{\partial (hv_{\theta,m}\alpha_i)}{\partial h}=-v_{\theta,m}\alpha_i, \quad\frac{\partial (hv_{\theta,m}\alpha_i)}{\partial (hv_{\theta,m})}=\alpha_i, \quad\frac{\partial (hv_{\theta,m}\alpha_i)}{\partial (h\alpha_l)}=v_{\theta,m}\delta_{i,l}.
\end{equation*}
Introducing the notation
\begin{equation*}
    V_{\boldsymbol{h,\alpha}}:=(h,hv_{r,m},h\alpha_1,\cdots,h\alpha_N),\quad V_{\boldsymbol{\gamma}}:=(hv_{\theta,m},h\gamma_1,\cdots,h\alpha_N),
\end{equation*}
with \( \boldsymbol{\alpha}=(\alpha_1,\cdots,\alpha_N)\) and with \( \boldsymbol{\gamma}=(\gamma_1,\cdots,\gamma_N)\), we obtain
\begin{equation*}
    \frac{\partial (hv_{r,m}\boldsymbol{\gamma}+hv_{\theta,m}\boldsymbol{\alpha})}{\partial V_{\boldsymbol{h,\alpha}}}\quad\myeq\quad \begin{pmatrix}
        -v_{r,m}\gamma_1-v_{\theta,m}\alpha_1 & \gamma_1 & v_{\theta,m} & & \\[6pt]
        & & & \ddots &  \\[6pt]
        & & & & v_{\theta,m} 
    \end{pmatrix},
\end{equation*} 
\begin{equation*}
    \frac{\partial (hv_{r,m}\boldsymbol{\gamma}+hv_{\theta,m}\boldsymbol{\alpha})}{\partial V_{\boldsymbol{\gamma}}}\quad\myeq\quad \begin{pmatrix}
        \alpha_1 & v_{r,m} & & \\[6pt]
        & & \ddots &  \\[6pt]
        & & & v_{r,m} 
    \end{pmatrix}.    
\end{equation*}

We denote the last term of the right-hand side of Equation \eqref{consflux} with the triple Legendre integral as \(\Tilde{F_i}:=h\sum_{j,k=1}^{N}A_{ijk}\alpha_j\alpha_k\). In \cite{HSWME}, it is shown that if \(j=1\), the only cases that need to be considered are \(k=i-1\) (if \(i\geq 2\)) and \(k=i+1\). Analogously, if \(k=1\), the only cases we need to consider are \(j=i-1\) (if \(i\geq 2\)) and \(j=i+1\). We can again consider four different situations.

\textbf{4.1.1 Case 1: \(i=1\).}
If \(i=1\), there are two cases for the index triplet \((i,j,k)\), with \(j=1\) or \(k=1\), which lead to a non-zero term: $(1,1,2) \text{ and } (1,2,1)$.

So we can pull these cases outside of the sum and we get:
\begin{equation*}
    \Tilde{F_1}=\sum_{j,k=1}^{N}A_{ijk}\alpha_j\gamma_k=A_{112}h\alpha_1\gamma_2+A{121}h\alpha_2\gamma_1+h\sum_{j,k=2}^{N}A_{ijk}\alpha_j\gamma_k.
\end{equation*}
The derivative of the last term with respect to $h$ reads
\begin{equation*}
    \frac{\partial \left(h\sum_{j,k=2}^{N}A_{ijk}\alpha_j\alpha_k\right)}{\partial h}=-\sum_{j,k=2}^{N}A_{ijk}\alpha_j\gamma_k.
\end{equation*}
After setting the higher order moments to zero, this term reduces to zero. The same can be observed when taking the derivative with respect to \(hv_{r,m}\), \(h\alpha_l\) and \(h\gamma_l\) (\(l=1,\ldots,N\)) and \(hv_{\theta,m}\). 

The partial derivatives of the first and second term denoted as \(\underline{\Tilde{F_1}}:=A_{112}h\alpha_1\gamma_2+A_{121}h\alpha_2\gamma_1\) are
\begin{align*}
    \frac{\partial\underline{\Tilde{F_1}}}{\partial(h\alpha_l)}&=A_{112}\gamma_2 \delta_{1,l}+A_{112}\gamma_1 \delta_{2,l}\quad \myeq \quad A_{112}\gamma_1 \delta_{2,l},\\[6pt]
    \frac{\partial\underline{\Tilde{F_1}}}{\partial(h\gamma_l)}&=A_{112}\alpha_1 \delta_{2,l}+A_{112}\alpha_2 \delta_{1,l}\quad \myeq \quad A_{112}\alpha_1 \delta_{2,l}.
\end{align*}
We denote the resulting entries by \(_{\alpha}c_1^F\) and \(_{\gamma}c_1^F\), respectively. %The first entry corresponds to the derivative with respect to \(h\alpha_2\) and the second entry corresponds to the derivative with respect to \(h\gamma_2\). 
%Note that the subscript \(\alpha\) refers to the fact that this entry has to with a derivative with respect to \(h\alpha_l\). Same for the subscript \(\gamma\). 
All other derivatives are zero, except for the derivative with respect to \(h\), but we can easily see that this derivative reduces to zero when we set higher order moments to zero.

\textbf{4.1.2 Case 2: \(i=2\)}.
Now, there are three cases for the index triplet \((i,j,k)\), with \(j=1\) or \(k=1\), which lead to a non-zero term: $(2,1,1), (2,1,3), \text{ and } (2,3,1)$.

We can pull this terms outside of the sum:
\begin{equation}
    \Tilde{F_2}=\sum_{j,k=1}^{N}A_{ijk}\alpha_j\gamma_k
    =A_{211}h\alpha_1\gamma_1+A_{213}h\alpha_1\gamma_3+A_{231}h\alpha_3\gamma_1+h\sum_{j,k=2}^{N}A_{ijk}\alpha_j\gamma_k.\label{consflux2}
\end{equation}
%\begin{align}
%    \Tilde{F_2}&=\sum_{j,k=1}^{N}A_{ijk}\alpha_j\gamma_k\\
%    &=A_{211}h\alpha_1\gamma_1+A_{213}h\alpha_1\gamma_3+A_{231}h\alpha_3\gamma_1+h\sum_{j,k=2}^{N}A_{ijk}\alpha_j\gamma_k.\label{consflux2}
%\end{align}
Again, the partial derivatives of the last term in the right-hand side of Equation \eqref{consflux2} all reduce to zero when setting the higher order moments to zero. Denoting the first three terms by \(\underline{\Tilde{F}_2}:=A_{211}h\alpha_1\gamma_1+A_{213}h\alpha_1\gamma_3+A_{231}h\alpha_3\gamma_1\), we obtain
\begin{equation*}
    \frac{\partial\underline{\Tilde{F}_2}}{\partial (h\alpha_l)}\quad \myeq \quad A_{211}\gamma_1\delta_{1,l}+A_{213}\gamma_1\delta_{1,3},\quad \frac{\partial\underline{\Tilde{F}_2}}{\partial (h\gamma_l)}\quad \myeq \quad A_{211}\alpha_1\delta_{1,l}+A_{213}\alpha_1\delta_{3,l}.
\end{equation*}
Considering the left equation, the first term leads to an entry \(_{\alpha}a_2^F\) and the second term leads to an entry \(_{\alpha}c_2^F\). Considering the right equation, the first term leads to an entry \(_{\gamma}a_2^F\) and the second term leads to an entry \(_{\gamma}c_2^F\).
Regarding the other derivatives, one can easily see that
\begin{equation*}
    \frac{\partial\underline{\Tilde{F}_2}}{\partial h}\quad \myeq \quad -A_{211}\alpha_1\gamma_1-\alpha_1\gamma_3,~ \frac{\partial\underline{\Tilde{F}_2}}{\partial (hv_{r,m})}=0,~ \frac{\partial\underline{\Tilde{F}_2}}{\partial (hv_{\theta,m})}=0.
\end{equation*}

\textbf{4.1.3 Case 3: \(2<i\leq N-1\)}.
In this case, there are four possibilities for the index triplet \((i,j,k)\): $(i,1,i-1), (i,1,i+1),  (i,i-1,1), \text{ and } (i,i+1,1)$.

We can pull this terms outside of the sum:
\begin{multline}
    \Tilde{F_i}=\sum_{j,k=1}^{N}A_{ijk}\alpha_j\gamma_k\\
    =A_{i,1,i-1}h\alpha_1\gamma_{i-1}+A_{i,1,i+1}h\alpha_1\gamma_{i+1}+A_{i,i-1,1}h\alpha_{i-1}\gamma_1 +A_{i,i+1,1}h\alpha_{i+1}\gamma_1+h\sum_{j,k=2}^{N}A_{ijk}\alpha_j\gamma_k. \label{consflux3}
\end{multline}
The derivatives of the last term of the right-hand side of Equation \eqref{consflux3} all reduce to zero when setting the higher order moments to zero. Denoting the other terms by \(\underline{\Tilde{F}_i}:=A_{i,1,i-1}h\alpha_1\gamma_{i-1}+A_{i,1,i+1}h\alpha_1\gamma_{i+1}+A_{i,i-1,1}h\alpha_{i-1}\gamma_1+A_{i,i+1,1}h\alpha_{i+1}\gamma_1\), we obtain
\begin{align*}
    \frac{\partial\underline{\Tilde{F}_i}}{\partial (h\alpha_l)}\quad &\myeq \quad A_{i,1,i-1}\gamma_1\delta_{i-1,l}+A_{i,1,i+1}\gamma_{i+1}\delta_{i+1,l},\\[6pt]
    \frac{\partial\underline{\Tilde{F}_i}}{\partial (h\gamma_l)}\quad &\myeq \quad A_{i,1,i-1}\alpha_1\delta_{i-1,l}+A_{i,1,i+1}\alpha_{i+1}\delta_{i+1,l}.
\end{align*}
Considering the first equation, the first term leads to an entry \(_{\alpha}a_i^F\), while the second term leads to en entry \(_{\alpha}c_i^F\). Analogously, the first term of the second equation leads to an entry \(_{\gamma}a_i^F\) and the second term of the second equation leads to an entry \(_{\gamma}a_i^F\). 

To summarize, the \(_{\alpha}a_i^F\)'s and the \(_{\alpha}c_i^F\) correspond to derivatives with respect to \(h\alpha_{i-1}\) and with respect to \(h\alpha_{i+1}\) of the conservative flux in the equation for \(\gamma_i\), respectively. The \(_{\gamma}a_i^F\)'s and the \(_{\gamma}c_i^F\) correspond to derivatives with respect to \(h\gamma_{i-1}\) and with respect to \(h\gamma_{i+1}\) of the conservative flux in the equation for \(\gamma_i\), respectively.

Furthermore, we find:
\begin{equation*}
    \frac{\partial\underline{\Tilde{F}_i}}{\partial h}\quad \myeq \quad 0,\quad\frac{\partial\underline{\Tilde{F}_i}}{\partial (hv_{r,m})}\quad \myeq \quad0,\quad\frac{\partial\underline{\Tilde{F}_i}}{\partial (hv_{\theta,m})}\quad \myeq \quad0.\\
\end{equation*}
%Thus, these derivatives do not lead to a non-zero entry in the flux Jacobian.

\textbf{4.1.4 Case 4: \(i=N\)}.
In this case, there are two possibilities for the index triplet \((i,j,k)\): $(N,1,N-1) \text{ and } (N,N-1,1)$.

Analogously to the previous cases, it can be easily shown that this case only leads to entries \(_{\alpha}a_N^F\), \(_{\alpha}c_N^F\), \(_{\gamma}a_N^F\) and \(_{\gamma}c_N^F\).

Using orthogonality and recursion formulas, we have:
\begin{align*}
    &_{\alpha}a_i^F=\frac{i}{2i-1}\gamma_1 \qquad \text{and} \qquad  _{\gamma}a_i^F=\frac{i}{2i-1}\alpha_1,  \qquad i=2,\ldots,N,\\[8pt]
    &_{\alpha}c_i^F=\frac{i+1}{2i+3}\gamma_1 \qquad \text{and} \qquad _{\gamma}c_i^F=\frac{i+1}{2i+3}\alpha_1, \qquad i=1,\ldots,N-1.
\end{align*}
In conclusion, we have for the modified conservative coefficient matrix:
\begin{equation*}
    \frac{\partial F}{\partial V_{\alpha}}_{mod}=
    \begin{pmatrix}
        & 1 & & & & & \\[6pt]
        gh-v_{r,m}^2-\frac{1}{3}\alpha_1^2 & 2v_{r,m} & \frac{2}{3}\alpha_1 & & & \\[6pt]
        -2 v_{r,m}\alpha_1 & 2\alpha_1 & 2v_{r,m} & c_1^F & & & \\[6pt]
        -\frac{2}{3}\alpha_1^2 & & a_2^F & 2v_{r,m} & \ddots & & \\[6pt]
        & & & \ddots & \ddots & c_{N-1}^F & \\[6pt]
        & & & & a_N^F & 2v_{r,m} \\[6pt]
        -v_{r,m}v_{\theta,m}-\frac{\alpha_1\gamma_1}{3} & v_{\theta,m} & \frac{\gamma_1}{3} & & & \\[6pt]
        -v_{r,m}\gamma_1-v_{\theta,m}\alpha_1 & \gamma_1 & v_{\theta,m} & _{\alpha}c_1^F & & \\[6pt]
        -\frac{2}{3}& &  _{\alpha}a_2^F & v_{\theta,m} & \ddots & \\[6pt]
        & & & \ddots & \ddots &  _{\alpha}c_{N-1}^F \\[6pt]
        & & & & _{\alpha}a_N^F & v_{\theta,m}
    \end{pmatrix},
\end{equation*}
and
\begin{equation*}
    \frac{\partial F}{\partial V_{\gamma}}_{mod}=
    \begin{pmatrix}
        & & & & & \\[6pt]
        & & & & & \\[6pt]
        v_{r,m} & \frac{\alpha_1}{3} & & & \\[6pt]
        \alpha_1 & v_{r,m} & _{\gamma}c_1^F & \\[6pt]
        &  _{\gamma}a_2^F & v_{r,m} & \ddots & \\[6pt]
        & & \ddots & \ddots &  _{\gamma}c_{N-1}^F \\[6pt]
        & & & _{\gamma}a_N^F & v_{r,m}
    \end{pmatrix},
\end{equation*}

where the first \(N+2\) rows of \(\frac{\partial F}{\partial V_{\gamma}}_{mod}\) are zero rows.

\textbf{4.2 The non-conservative part \(\boldsymbol{Q_i}\)}.
The equations for \(h,hv_{r,m}\) and \(hv_{\theta,m}\) do not contain non-conservative terms, so the non-conservative system matrix will contain zero rows for these equations. As observed above, the non-conservative part in the equations for \(h\alpha_i\), \(i=1,\ldots,N\), does not depend on derivatives with respect to \(h\gamma_l\), \(l=1,\ldots,N\). 
Recall the non-conservative part in the equation for \(h\gamma_i\):
\begin{equation}\label{nonConservativeFlux}
    Q_i=v_{\theta,m}\frac{\partial(h\alpha_i)}{\partial r}-\sum_{j,k=1}^{N}B_{ijk}\gamma_k\frac{\partial(h\alpha_j)}{\partial r}.
\end{equation}
Thus, the equations for \(h\gamma_i\), \(i=1,\ldots,N\), do not depend on derivatives with respect to \(h\gamma_l\), \(l=1,\ldots,N\), either. From this, it follows that we can write the modified non-conservative terms (i.e., with higher order moments set to zero) as 
\begin{equation*}
    Q_{mod}=
    \begin{pmatrix}
        Q_1 & Q_2 \\
        Q_3 & Q_4
    \end{pmatrix},
\end{equation*}
in which \(Q_2\in\mathbb{R}^{(N+3)\times N}\) and \(Q_4\in\mathbb{R}^{N\times N}\) are zero matrices. \(Q_1\in\mathbb{R}^{(N+3)\times (N+3)}\) contains the terms corresponding to the non-conservative flux for the first \(N+2\) equations with respect to the first \(N+3\) partial derivatives, see Equation \eqref{nonconsflux}, and a zero row corresponding to the equation for \(hv_{\theta,m}\). The remaining entries to compute are the entries of the matrix \(Q_3 \in\mathbb{R}^{N\times (N+3)}\). This matrix corresponds to the equation for \(\gamma_i\), \(i=1,\ldots,N\), and to the derivatives with respect to \(h,hv_{r,m}\) and \(h\alpha_l\) (\(l=1,\ldots,N)\). 
The first term of the non-conservative flux leads to the term \(v_{\theta,m}\) on the diagonal of \(Q_3\). The second term is simplified when we set higher order moments to zero:
\begin{equation*}
    \sum_{j=1}^{N}(2i+1)\left(\int_0^1\phi_i'\left( \int_0^1\phi_jd\hat{\zeta}\right)\phi_1d\zeta\right)\gamma_1\frac{\partial(h\alpha_j)}{\partial r}.
\end{equation*}
As stated before, it can be proved that \(B_{ij1}=0\) for \(|i-j|\neq 1\), except for the first two rows \cite{HSWME}. \(Q_i\) does not depend on any derivative with respect to \(h\), \(hv_{r,m}\) and \(hv_{\theta,m}\), so the entries of the corresponding columns will be zero. We can see that the second term in \eqref{nonConservativeFlux} only leads to non-zero entries \(a_i^Q\) on the first lower diagonal and non-zero entries \(c_i^Q\) on the first upper diagonal:
\begin{align*}
    &a_i^Q=\frac{i+1}{2i-1}\gamma_1, \qquad i=2,\ldots,N,\\[8pt]
    &c_i^Q=\frac{i}{2i+3}\gamma_1, \qquad i=1,\ldots,N-1.
\end{align*}

So \(Q_3\) has the following form:
\begin{equation*}
    Q_3=
    \begin{pmatrix}
        & & & v_{r,m} & c_1^Q & & \\[6pt]
        & & & a_2^Q & v_{r,m} & \ddots & & \\[6pt]
        & & & & \ddots & \ddots & c_{N_r-1}^Q & \\[6pt]
        & & & & & a_{N_r}^Q & v_{r,m} &\\[6pt]
    \end{pmatrix}.
\end{equation*}

From this, the matrix \(Q_{mod}\) can be constructed and the modified system matrix is 
\begin{equation*}
    A_{HA}^{(N,N)}=\frac{\partial_F}{\partial V}_{mod}-Q_{mod}.
\end{equation*}
This completes the proof.
\end{proof}

%\clearpage

\section{System Matrix Of Second Order Axisymmetric System With Full Velocity Expanded}
\label{app:B}
\begin{comment}
The system matrix of the \((2,2)\)th order ASWME system is given by:
\begin{equation*}
    A_A^{(2,2)}=
    \begin{pmatrix}
        0 & 1 & 0 & 0 & 0 & 0 & 0 \\
        d_1 & 2 v_{r,m} & \frac{2 \alpha _1}{3} & \frac{2 \alpha _2}{5} & 0 & 0 & 0 \\
        d_2 & 2 \alpha _1 & v_{r,m}+\alpha _2 & \frac{3 \alpha _1}{5} & 0 & 0 & 0 \\
        d_3 & 2 \alpha _2 & \frac{\alpha _1}{3} & v_{r,m}+\frac{3 \alpha _2}{7} & 0 & 0 & 0 \\
        \Tilde{d}_1 & v_{\theta ,m} & \frac{\gamma _1}{3} & \frac{\gamma _2}{5} & v_{r,m} & \frac{\alpha _1}{3} & \frac{\alpha _2}{5} \\
        \Tilde{d}_2 & \gamma _1 & \frac{3 \gamma _2}{5} & \frac{\gamma _1}{5} & \alpha _1 & v_{r,m}+\frac{2 \alpha _2}{5} & \frac{2\alpha_1}{5} \\
        \Tilde{d}_3 & \gamma _2 & -\frac{\gamma _1}{3} & \frac{\gamma _2}{7} & \alpha _2 & \frac{2 \alpha _1}{3} & v_{r,m}+\frac{2 \alpha_2}{7}
    \end{pmatrix},
\end{equation*}
where
\end{comment}

The entries in the first column of the system matrix of the \((2,2)\)th order ASWME system in Equation \eqref{A22} are given by:
\begin{align*}
    d_1 &= gh-v_{r,m}^2-\frac{\alpha_1^2}{3}-\frac{\alpha_2^2}{5},\\[5pt] 
    d_2 &= -v_{r,m}-\frac{4}{5}\alpha_1\alpha_2,\\[5pt] 
    d_3 &= -\frac{2}{21} \left(3 \alpha _2 \left(7 v_{r,m}+\alpha _2\right)+7 \alpha _1^2\right),\\[5pt]
    \Tilde{d}_1 &= -v_{r,m} v_{\theta ,m}-\frac{1}{3} \alpha _1 \gamma _1-\frac{\alpha _2 \gamma _2}{5} ,\\[5pt] 
    \Tilde{d}_2 &= \frac{1}{5} \left(-\gamma _1 \left(5 v_{r,m}+2 \alpha _2\right)-\alpha _1 \left(5 v_{\theta ,m}+2
   \gamma _2\right)\right),\\[5pt] 
    \Tilde{d}_3 &= -\gamma _2 v_{r,m}-\frac{1}{7} \alpha _2 \left(7 v_{\theta ,m}+2 \gamma _2\right)-\frac{2}{3} \alpha _1
   \gamma _1.
\end{align*}